\documentclass[a4paper,10pt]{article}

\usepackage{graphicx}
\usepackage{amsfonts}
\usepackage{amstext}
\usepackage{amssymb,amsmath,amsbsy,amsthm}
\usepackage[T1]{fontenc}
\usepackage[utf8]{inputenc}
\usepackage{authblk}
\usepackage{caption}
\usepackage{subcaption}
\usepackage{bm}
\usepackage{hyperref}

\newtheorem{theorem}{Theorem}

\newcommand{\R}{\mathbb{R}}
\newcommand{\0}{\mathbf{0}}
\newcommand{\1}{\mathbf{1}}
\newcommand{\cc}{\mathbf{c}}
\newcommand{\ddd}{\mathbf{d}}
\newcommand{\q}{\mathbf{q}}
\newcommand{\rrr}{\mathbf{r}}
\newcommand{\x}{\mathbf{x}}
\newcommand{\X}{\bm{X}}
\newcommand{\uuu}{\mathbf{u}}
\newcommand{\vvv}{\mathbf{v}}

\newcommand{\y}{\mathbf{y}}

\newcommand{\z}{\mathbf{z}}

\newcommand{\I}{\bm{I}}
\newcommand{\LL}{\bm{L}}
\newcommand{\M}{\bm{M}}

\newcommand{\W}{\bm{W}}

\newcommand{\C}{\bm{C}}
\newcommand{\DD}{\bm{D}}

\newcommand{\tr}{\mathrm{tr}}
\newcommand{\diag}{\mathrm{diag}}

\newcommand{\rk}{\mathrm{rank}}

\newcommand{\Vol}{{\mathrm{Vol}}}

\newcommand{\md}{\mathrm{md}}

\title{Regularity based spectral clustering and mapping the Fiedler-carpet} 

\author{Marianna Bolla \thanks{Department of Stochastics (DS), Budapest University of Technology and Economics (BME). e-mail: marib@math.bme.hu} 	\and Vilas Winstein  \thanks{Renyi Institute of Mathematics, Budapest, Hungary, e-mail: vilas@renyi.hu} \and Tao You  \thanks{Department of Mathematics, Middlebury College, e-mail: tyou@middlebury.edu} \and Frank Seidl  \thanks{Department of Mathematics, University of Michigan, e-mail: fcseidl@umich.edu} \and Fatma Abdelkhalek \thanks{DS, BME and Faculty of Commerce, Assiut University, Egypt, e-mail: fatma@math.bme.hu}
}

\begin{document}
\maketitle

\section*{Abstract}

Spectral clustering is discussed from many perspectives, by extending it to
rectangular arrays
and discrepancy minimization too. Near optimal clusters are obtained
with singular value decomposition and with the weighted $k$-means algorithm.
In case of rectangular arrays, this means
enhancing the method of correspondence analysis with clustering, 
and in case of  edge-weighted graphs,
a normalized Laplacian based clustering. In the latter case it is proved
that a spectral gap between the $(k-1)$th and $k$th smallest positive
eigenvalues of the normalized Laplacian
matrix gives rise to a sudden decrease of the inner cluster variances when
the number of clusters of the vertex representatives is $2^{k-1}$, 
but only the first $k-1$ eigenvectors, constituting the so-called
Fiedler-carpet, are used in the representation.
Application to directed migration graphs is also discussed.

\noindent
\textbf{Keywords}: multiway discrepancy, correspondence analysis, normalized
Laplacian, multiway cuts, Fiedler-vector, weighted $k$-means algorithm.

\noindent
\textit{MSC2020}: 05C50,62H25,65H10.

\section{Introduction}

Spectral graph theory started developing about 50 years ago 
(see, e.g. A.J. Hoffman~\cite{Hoffman2}, M.~Fiedler~\cite{Fiedler},
D.M. Cvetkovi\'c~\cite{Cvetkovic}, and F.~Chung~\cite{Chung0}) to characterize 
certain structural properties 
of a graph by means of the eigenvalues of its adjacency or Laplacian matrix.
Later on, in the last decades of the 20th century, 
the eigenvectors corresponding to some near zero eigenvalues of
the Laplacian matrix were also used for clustering  the vertices into
disjoint parts so that the inter-cluster relations are negligible compared
to the intra-cluster ones (usual purpose of  cluster analysis as a machine
learning technique). 
In this setup, the famous \textit{Fiedler-vector}, the  eigenvector, 
corresponding to the smallest positive Laplacian eigenvalue was used to
classify the vertices into two parts; while later, $k$ eigenvectors
($k\ge 2$) entered into the $k$-clustering task. 
It was also noted in~\cite{Alpert} that the more eigenvectors, the better,
but not exact estimates between the spectral gap and the quality of
the clustering were available, except the $k=1$ and $k=2$ cases; former
related to the isoperimetric number and expander graphs
(see e.g.~\cite{Mohar1,Chung0,Hoory}) and latter to the sum of the inner 
variances of 2-clusterings (see~\cite{Bolla3}).

Since then, the problem was generalized in several ways:
to edge-weighted graphs and rectangular arrays of nonnegative entries (e.g.
microarrays in biological genetics and 
forensic science~\cite{Bolla10,Bolla16,BollobasN,Kluger}),
and to degree-corrected adjacency and Laplacian 
matrices~\cite{Bolla13,Chung2}. 
After the millennium, 
physicists and social scientists introduced modularity matrices
and investigated so-called anti-community structures (intra-cluster relations 
are negligible compared to the inter-cluster ones) in contrast to the former
community structures~\cite{Newman10}. 
By uniting these two approaches,  so-called regular cluster pairs with
small discrepancy can be defined, where homogeneous clusters are looked for
(e.g. in microarrays one looks for groups of genes that similarly influence
the same groups of conditions)~\cite{Bolla16}. 
The existence of such a regular structure is
theoretically guaranteed by the Abel-prize winner 
Szemeredi's regularity lemma~\cite{Szemeredi}
that for any small positive $\varepsilon$  guarantees a universal number $k$ 
of clusters (irrespective of the numbers 
of vertices) such that partitioning the vertices into $k$ parts (and
possibly a small
exceptional one), the pairs have discrepancy less than $\varepsilon$.  
However, this $k$ can be enormously large and not applicable to
practical purposes. Our purpose is to give a moderate $k$
where the sum of the inner variances of
$2^{k-1}$ clusters is estimated from above by the spectral gap between the
$(k-1)$th and $k$th positive normalized Laplacian eigenvalues
in their non-decreasing order,
even in the worst case scenario. This is the generalization of a theorem
of~\cite{Bolla3} that was applicable only to the $k=2$ situation.
There $2^{k-1}$ and $k$ were the same, and the Fiedler-vector was used for
clustering into $k=2$ parts. If $k>2$, then we use $k-1$ non-trivial 
eigenvectors that form what we call \textit{Fiedler-carpet}.
In special cases, e.g. in case of generalized multiclass random or quasirandom
graphs,  both the objective function of the $k$-means algorithm and the
the $k$-way discrepancy dramatically decreases 
compared to the $k-1$ one~\cite{Bolla20}, where the number of clusters is 
one more than the number of eigenvectors used in the representation.
However, in the generic case, the number of eigenvectors entered into the 
classification is much less than the number of clusters,
 that is also a good news
from the point of view of computational complexity.

The organization of the paper is as follows.
In Section~\ref{pre}, the most important notions and facts are defined,
concerning normalized Laplacian spectra and multiway discrepancy of
rectangular arrays with nonnegative 
entries~\cite{Alon10,Bolla13,Bolla16,Bolla20,BollobasN}.
In Section~\ref{main}, the main theorem is stated and proved.
The proof is based on the energy minimizing representation and
 analyzing the structure of the vertex representatives
and the underlying spectral subspaces that are mapped in a convenient way.
In Appendix A, next to technical considerations, 
simulation results, and plots of the so-obtained Fiedler-carpet 
are also presented.
Real-life application to the directed graph of migration data is discussed in 
Section~\ref{appl}, supplemented with more figures in Appendix B.
In Section~\ref{sum}, the main contributions are summarized and
conclusions are drawn.

\section{Minimal and regular cuts versus spectra}\label{pre}

Let $\W$ be the $n\times n$ edge-weight matrix of a graph $G$ on $n$ vertices.
It is symmetric, has 0 diagonal and nonnegative entries. In case of a
simple graph, $\W=(w_{ij})$ is the usual adjacency matrix. 
The generalized degrees are $d_i = \sum_{j=1}^n w_{ij}$, for $i=1,\dots ,n$.
Assume that $d_i$s are all positive and the diagonal degree-matrix $\DD$ 
contains  them in its main diagonal. The \textit{Laplacian} of $G$ is
$\LL =\DD -\W$, while its \textit{normalized Laplacian} is  
$$
 \LL_{\DD } =\I -\DD^{-1/2} \W \DD^{-1/2} .
$$
Because of the normalization, $\LL_{\DD}$ is not affected by the scaling of the
edge-weights, therefore $\sum_{i=1}^n d_i =1$ can be assumed.   
$\LL_{\DD}$  is positive semidefinite, and if $\W$ is irreducible 
($G$ is connected), then its eigenvalues are
$$
  0=\lambda_0 <\lambda_1  \le\dots \le \lambda_{n-1} \le 2
$$ 
with unit-norm pairwise orthogonal eigenvectors $\uuu_0 ,\uuu_1 ,\dots ,\uuu_{n-1}$.
In particular, $\uuu_0 = (\sqrt{d_1} ,\dots ,\sqrt{d_n } )^T =: \sqrt{\ddd }^T$.

For $1<d <n$, the row vectors of the $n\times d$ matrix $\X^* =(\DD^{-1/2} \uuu_1 ,
\dots , \DD^{-1/2} \uuu_d )$
are optimal $d$-dimensional representatives $\rrr^*_1 ,\dots ,\rrr^*_n$ of the 
vertices that minimize the \textit{energy} function
\begin{equation}\label{tr}
 Q_d (\X ) = \sum_{i=1}^{n-1} \sum_{j=i+1}^n  w_{ij} \| \rrr_i -\rrr_j \|^2 =
  \tr (\X^T \LL \X ) ,
\end{equation}
where the general vertex-representatives $\rrr_1 , \dots ,\rrr_n \in \R^d$
are row vectors of the $n\times d$  matrix  $\DD^{-1/2} \X$. They satisfy
the constraints $\X^T \DD \X =\I_d$ and $\sum_{i=1}^n d_i \rrr_i =\0$.
Minimizing the energy $Q_d (\X )$ supports representatives of vertices
connected with large weights to be close to each other.  
The minimum of $Q_d (\X )$  is the sum of the
$d$ smallest positive eigenvalues of $\LL_{\DD}$.

 The \textit{weighted k-variance} of these representatives is defined as
\begin{equation}\label{wkszoras}
 {S}_k^2  (\X ) =\min_{(V_1 ,\dots ,V_k ) \in {\cal P}_k }
\sum_{i=1}^k \sum_{j\in V_i } d_j \| \rrr_j -{ \cc }_i \|^2 ,
\end{equation}
where $\Vol (U ) =\sum_{j\in U} d_j$,  
${\cc }_i =\frac1{\Vol (V_i ) } \sum_{j\in V_i } d_j \rrr_j $ is the
weighted center of the cluster $V_i$, and the minimization is over 
proper $k$-partitions $P_k =(V_ 1 ,\dots ,V_k )$ of the vertex-set, their
collection is denoted by ${\cal P}_k$.

It is the weighted $k$-means algorithm that approaches this minimum.
In~\cite{Ostrovsky}, it is stated that if the data
satisfy the $k$-clusterable criterion ($S_k^2 \le \epsilon^2 S_{k-1}^2$
with a small enough $\epsilon$), then there is a PTAS (polynomial time
approximation scheme) for the $k$-means problem.
This is the situation we usually owe.

It is well known that the $k$ bottom eigenvalues of the normalized Laplacian 
matrix estimate the $k$-way \textit{normalized cut} of $G$ which is
$f_k (G) =\min_{P_k \in {\cal P}_k } f (P_k ,G )$, where
$$
\begin{aligned}
 f (P_k , G )&= \sum_{a=1}^{k-1} \sum_{b=a+1}^k \left( \frac1{ \Vol (V_a) }+
   \frac1{ \Vol (V_b ) } \right) w (V_a ,V_b ) \\
 &=\sum_{a=1}^{k} \frac{w (V_a ,{\overline V}_a )}{ \Vol (V_a ) } 
  =k-\sum_{a=1}^{k} \frac{w (V_a ,V_a )}{ \Vol (V_a ) } .
\end{aligned}
$$
Here $ w (V_a ,V_b ) = \sum_{i \in V_a} \sum_{j \in V_b} w_{ij}$ is the
weighted cut between the cluster pairs.
Since $\sum_{i=1}^{k-1} \lambda_i$ 
is the overall minimum of $Q_k$ (on the orthogonality
constraints) and $f_k (G)$ is the minimum over partition vectors (having
stepwise constant coordinates over the parts of $P_k$), the relation
\begin{equation}\label{normcut}
 \sum_{i=1}^{k-1} \lambda_i \le f_k (G)
\end{equation}
is easy to prove. This estimate is sharper if the eigensubspace
spanned by the corresponding eigenvectors is closer to that of the partition
vectors in the convenient $k$-partition of the vertices, the one produced
by the \textit{weighted k-means algorithm}. Therefore $S_k^2 (\X_{k-1}^* )$ 
indicates the
quality of the $k$-clustering based on the $k-1$ bottom eigenvectors
(except the trivial one). 
Later it will be used that neither $Q_k$ nor  $S_k^2 (\X_{k-1}^* )$ is affected
by the orientation of the orthonormal eigenvectors.  

In~\cite{Bolla3} it is proved that
$S_2^2 (\X_1^* ) \le \frac{\lambda_1}{\lambda_2}$, 
so the larger the gap after the
first positive eigenvalue of $\LL_{\DD}$, the sharper the estimate
in~\eqref{normcut} is. Here this statement is generalized  to the gap
between $\lambda_{k-1}$ and $\lambda_k$, but $2^{k-1}$
clusters are considered based on $(k-1)$-dimensional vertex representatives. 
In the literature (see,
e.g.~\cite{Luxburg,Ng,Shi}) the number of clusters is usually the 
same as the number of eigenvectors entered into the classification. 
The message of Theorem~\ref{gap} of Section~\ref{main} is
that the number of clusters is much higher in the generic case than
the dimension of the representatives, at least in the minimum multiway
cut problems. Though, with discrepancy objective, the famous
Szemer\'edi Regularity Lemma~\cite{Szemeredi} also suggests this.

Now consider the discrepancy view and the related matrices.
The whole can better be illustrate on rectangular arrays of
nonnegative entries; simple, edge-weighted, and directed graphs are
special cases.

In many applications, for example when microarrays are analyzed, our
data are collected in the form of an $m\times n$ rectangular matrix
$\C=(c_{ij})$ of
nonnegative real entries. (If the entries are integer frequency counts, then
the array is called contingency table in statistics.) 
We assume that $\C$ is
\textit{non-degenerate}, i.e.
$\C \C^T$ (when $m\le n$) or $\C^T \C$  (when $m > n$) is \textit{irreducible}.
Consequently, the row-sums
$d_{row,i} =\sum_{j=1}^n c_{ij}$ and column-sums $d_{col,j}=\sum_{i=1}^m c_{ij}$
of $\C$ are strictly positive, and the diagonal matrices
$\DD_{row} =\diag (d_{row,1} ,\dots ,d_{row,m})$ and 
$\DD_{col} =\diag (d_{col,1} ,\dots ,d_{col,n})$ are regular.
Without loss of generality, we also assume that
$\sum_{i=1}^n \sum_{j=1}^m c_{ij} =1$, since neither our main object, the
\textit{normalized contingency table} 
\begin{equation}\label{cnor}
 \C_{\DD} =  \DD_{row}^{-1/2} \C  \DD_{col}^{-1/2}  ,
\end{equation}
nor the \textit{multiway discrepancies} to be introduced  
are affected by the scaling of the entries of $\C$.
It is known  that the singular values  of  $\C_{D}$ are in the [0,1]
interval. The positive ones, enumerated in non-increasing order, are the real 
numbers
$$
 1=s_0 >s_1 \ge \dots \ge s_{r-1} > 0  ,
$$
where $r= \rk (\C_D )=\rk (\C )$. Provided $\C$ is non-degenerate,
1 is a single singular value; 
it will be called trivial and denoted
by $s_0$, since it corresponds to the trivial singular vector pair,
which are disregarded in the clustering problems. 
This is a well-known fact of \textit{correspondence analysis}, for further
explanation see~\cite{Bolla13,Faust} and the subsequent paragraph.

For a given integer $1\le k\le \min \{ m,n \}$,
 we are looking for $k$-dimensional
representatives $\rrr_1 ,\dots ,\rrr_m \in \R^k$ of the rows and 
 $\q_1 ,\dots ,\q_n \in \R^k$ of the columns
such that they minimize the energy function
\begin{equation}\label{objrec}
 Q_k = \sum_{i=1}^m \sum_{j=1}^n  c_{ij} \| \rrr_i -\q_j \|^2 
\end{equation}
subject to
\begin{equation}\label{constrec}
 \sum_{i=1}^m d_{row,i} \rrr_i \rrr_i^T =\I_k  \quad \textrm{and} \quad
 \sum_{j=1}^n d_{col,j} \q_j \q_j^T =\I_k .
\end{equation}

When minimized, the objective function $Q_k$ favors $k$-dimensional placement
of the rows and columns such that representatives of columns and rows with 
large association are
forced to be close to each other. 
It is easy to prove that the minimum is obtained by the
singular value decomposition (SVD)
\begin{equation}\label{correspmatrix}
 \C_{\DD} =\sum_{k=0}^{r-1} s_k \vvv_k \uuu_k^T ,
\end{equation}
where $r\le \min \{ n,m \}$ is the rank of $\C_{\DD}$.
The constrained minimum of $Q_k$ is
$2k -\sum_{i=0}^{k-1} s_i$ and it is 
attained with row- and column-representatives that are row vectors of the
matrices
$\DD_{row}^{-1/2} (\vvv_0 ,\vvv_1 ,\dots ,\vvv_{k-1} )$ 
and $\DD_{col}^{-1/2} (\uuu_0 ,\uuu_1 ,\dots ,\uuu_{k-1} )$, respectively.

Note that if the entries of $\C$ are frequency counts and their  sum ($N$)
is the sample size, then
the $\chi^2$ statistic, which measures the deviation from \textit{independence},
is 
\begin{equation}\label{chi2}
  \chi^2 =N\sum_{i=1}^{r-1}s_i^2  .
\end{equation} 
If the $\chi^2$ test based on this statistic indicates significant deviance
from independence (i.e. from the rank 1 approximation of $\C$), 
then one may look for 
rank $k$ approximation $(1<k<r=\rk )$, which  is constructed by the first $k$
singular vector pairs.

With bi-clustering the rows and columns of $\C$ one also may look for sub-tables
close to independent ones. This is measured by the \textit{discrepancy}.
In~\cite{Bolla16},
the multiway discrepancy  of the rectangular matrix $\C$ of nonnegative entries
in the proper $k$-partition $R_1 ,\dots ,R_k$ of its rows and
$C_1 ,\dots ,C_k$ of its columns is defined as
\begin{equation}\label{disk} 
\begin{aligned}
 \md (\C ; R_1 ,\dots ,R_k , C_1 ,\dots ,C_k ) &=
 \max_{\substack{1\le a , b\le k \\X\subset R_a , \, Y\subset C_b}} 
 \frac{|c (X, Y)-\rho (R_a,C_b ) \Vol (X)\Vol (Y)|}{\sqrt{\Vol(X)\Vol(Y)}} \\
 &=|\rho (X, Y)-\rho (R_a,C_b ) | {\sqrt{\Vol(X)\Vol(Y)}}
\end{aligned}
\end{equation}
where
$c (X, Y) =\sum_{i\in X} \sum _{j\in Y} c_{ij}$ is the cut between
$X\subset R_a$ and $Y\subset C_b$, 
$\Vol (X) = \sum_{i\in X} d_{row,i}$ is the volume of the row-subset $X$, 
$\Vol (Y) = \sum_{j\in Y} d_{col,j}$ is the volume of the column-subset $Y$, 
whereas
$\rho (R_a,C_b) =\frac{c(R_a,C_b)}{ \Vol (R_a) \Vol (C_b)}$ denotes the relative
density between $R_a$ and $C_b$.
The minimum $k$-way discrepancy  of  $\C$ itself is
$$
 \md_k (\C ) = \min_{\substack{R_1 ,\dots ,R_k \\ C_1 ,\dots ,C_k } } 
 \md (\C ; R_1 ,\dots ,R_k , C_1 ,\dots ,C_k ).
$$
In~\cite{Bolla16} the following is proved: 
$$
 s_k  \le 9\md_{k } (\C )  (k+2 -9k\ln \md_{k } (\C )) ,
$$
provided $0<\md_{k } (\C ) <1$.
In the forward direction, the
following is established in~\cite{Bolla13}. Given the $m\times n$ contingency
table $\C$, consider the spectral clusters $R_1 ,\dots ,R_k$ of its rows and 
$C_1 ,\dots ,C_k$ of its columns,
obtained by applying the weighted $k$-means algorithm 
to the $(k-1)$-dimensional  row- and column representatives.
Let $S_{k,row}^2$ and 
$S_{k,col}^2$  denote the minima of the $k$-means algorithm with them, 
respectively.
Then, under some balancing conditions for the margins 
and for the cluster sizes,  
$\md_k (\C ) \le B (\sqrt{2k} (S_{k,row} +S_{k,col}) +s_k )$,
with some constant $B$, which depends only on the constants of the balancing 
conditions, and does not depend on $m$ and $n$.
Roughly speaking, the two directions
together imply that if $s_k$ is `small' and `much smaller' than $s_{k-1}$,
then one may
expect a simultaneous $k$-clustering of the rows and columns of $\C$
with small $k$-way discrepancy. This is the case of generalized random
and quasirandom graphs~\cite{Bolla20}.

This notion can be extended to an edge-weighted graph $G$
and denoted  by $\md_k (G)$. In that setup, $\C$ plays the role of 
the weighted adjacency matrix
(symmetric in the undirected; quadratic, but usually not symmetric 
in the directed case).
Here the singular values of the normalized adjacency matrix are the absolute 
values of the eigenvalues, which enter into the
estimates, in decreasing order.

At this point, some new matrices are introduced, originally defined 
by physicists  (see~\cite{Newman10}). The \textit{modularity matrix}
of an edge-weighted  graph $G$ is defined as $\M =\W-\ddd \ddd^T $,
where the entries of $\W$  sum to 1.
The \textit{normalized modularity matrix} of $G$ (see~\cite{Bolla11}) is
$$
  \M_{\DD} =\DD^{-1/2} \M \DD^{-1/2} = \DD^{-1/2}\W \DD^{-1/2}-\sqrt{\ddd }
 \sqrt{\ddd }^T =\W_{\DD}  -\sqrt{\ddd }\sqrt{\ddd }^T .
$$
The normalized modularity matrix is the normalized edge-weight matrix deprived
of the trivial dyad. Obviously, $\LL_{\DD } = \I -\W_{\DD} = \I -\M_{\DD} -
\sqrt{\ddd }\sqrt{\ddd }^T$, see~\cite{Bolla13}.

Therefore, the  $k-1$ largest singular values of $\M_{\DD }$ are
the absolute values of the $k-1$ largest eigenvalues of $\W_{\DD }$
(except the trivial 1).
These, in turn, are 1 minus the  $k-1$ positive eigenvalues  of 
$\LL_{\DD}$ which are
in the farthest distance from 1. If those are all less then 1, then these
are $1-\lambda_1 ,\dots , 1-\lambda_{k-1}$. In this case, the regularity
based spectral clustering boils down to the minimum cut objective.
 
Otherwise, for a $1<k<n$ fixed integer,
in the modularity based spectral clustering,
we look for the proper $k$-partition $V_1,\dots ,V_k $ of the vertices
such that the within- and between cluster discrepancies are minimized.
To motivate the introduction of the exact discrepancy measure observe that
the $ij$ entry of $\M$ is $w_{ij} - d_i d_j$, which is the difference between
the actual connection of the vertices $i,j$ and the connection that is
expected under independent attachment of them with probabilities $d_i$ and
$d_j$, respectively. Consequently, the difference between the actual
and the expected connectedness of the subsets $X,Y \subset V$ is
$$
 \sum_{i\in X} \sum_{j\in Y} (w_{ij} - d_i d_j ) =w(X,Y) -\Vol (X) \Vol (Y).
$$

A directed edge-weighted graph $G=(V,\W )$ is described by its quadratic, but 
usually not symmetric
weighted adjacency  matrix $\W=(w_{ij})$ of zero diagonal, 
where $w_{ij}$ is the nonnegative weight of the $j\to i$ edge $(i \ne j )$.
The row-sums
$d_{in,i} =\sum_{j=1}^n w_{ij}$ and column-sums $d_{out,j}=\sum_{i=1}^n w_{ij}$
of $\W$ are the \textit{in- and out-degrees}, while 
$\DD_{in} =\diag (d_{in,1} ,\dots ,d_{in,n})$ and 
$\DD_{out} =\diag (d_{out,1} ,\dots ,d_{out,n})$ are the diagonal
in- and out-degree matrices.
The multiway discrepancy  of the directed,
edge-weighted graph $G=(V,\W )$ in the in-clustering
$V_{in,1} ,\dots ,V_{in,k}$ and out-clustering $V_{out,1} ,\dots ,V_{out,k}$ 
of its vertices is
$$
\begin{aligned}
 &\md (G; V_{in,1} ,\dots ,V_{in,k}, V_{out,1} ,\dots ,V_{out,k})  \\
 &=\max_{\substack{1\le a , b\le k \\ X\subset V_{out,a} , \, Y\subset V_{in,b}}} 
 \frac{|w (X, Y)-\rho (V_{out,b},V_{in,a} ) \Vol_{in} (X)\Vol_{out} (Y)|} 
 {\sqrt{\Vol_{n}(X)\Vol_{out} (Y)}},
\end{aligned}
$$
where $w(X,Y)$ is the sum of the weights of the $X\to Y$ edges, whereas
$\Vol_{in} (X) =\sum_{i\in X} d_{in,i}$ and 
$\Vol_{out} (Y) =\sum_{j\in Y} d_{out,j}$ are the out- and in-volumes,
respectively. 
The minimum $k$-way discrepancy  of the directed 
edge-weighted graph $G=(V,\W )$ is
$$
 \md_k (G) = \min_{\substack{V_{in,1} ,\dots ,V_{in,k} \\ 
  V_{out,1} ,\dots ,V_{out,k} }}
 \md (G; V_{in,1} ,\dots ,V_{in,k}, V_{out,1} ,\dots ,V_{out,k}).
$$
In~\cite{Bolla16} it is proved that 
$$
 s_k  \le 9\md_{k } (G )  (k+2 -9k\ln \md_{k } (G )) ,
$$
where $s_k$ is the $k$-th largest nontrivial singular value  of
the normalized weighted adjacency matrix 
$\W_{\DD}=\DD_{in}^{-1/2} \W \DD_{out}^{-1/2}$.
In Section~\ref{appl} we apply the SVD of $\W_{\DD }$ to find migration patterns,
i.e. emigration and immigration trait clusters. 

\section{Mapping the Fiedler-carpet: more clusters than
eigenvectors}\label{main}

\begin{theorem}\label{gap}
Let $G=(V,\W )$ be connected edge-weighted graph 
with generalized degrees $d_1 ,\dots ,
d_n$ and assume that $\sum_{i=1}^n d_i =1$.
Let  $0=\lambda_0 <\lambda_1 
\le\dots \le \lambda_{n-1} \le 2$ denote the eigenvalues of the normalized 
Laplacian 
matrix $\LL_{\DD}$ of $G$. Then for the weighted $2^{k-1}$-variance of the 
optimal $(k-1)$-dimensional vertex
representatives, comprising row vectors of the matrix $\X_{k-1}^*$,
the following upper estimate holds:
$$
 {S}_{2^{k-1}}^2(\X_{k-1}^* ) \le \dfrac{\sum_{j=1}^{k-1}\lambda_j}{\lambda_k} ,
$$
provided $\lambda_{k-1} <\lambda_k$.
\end{theorem}

\noindent
\begin{proof}
Recall that, with the notation of Section~\ref{pre},  
$\X_{k-1}^* =(\DD^{-1/2} \uuu_1 , \dots , \DD^{-1/2} \uuu_{k-1})$,
where the trivial $\DD^{-1/2} \uuu_0 =\1$ vector is disregarded, and 
$\uuu_0 , \uuu_1 ,\dots \uuu_{k-1}$  are unit-norm, 
pairwise orthogonal eigenvectors
corresponding to the eigenvalues $0=\lambda_0 < \lambda_1 \le  \dots \le 
\lambda_{k-1}$ of $\LL_D$, respectively. As the trivial dimension is disregarded,
we only use the
coordinates of the vectors $\x_j := \DD^{-1/2} \uuu_j =
(x_{j1} ,\dots ,x_{jn})^T$ for $j=1,\dots ,k-1$. 

Since $\uuu_1 ,\dots \uuu_{k-1}$ form an orthonormal set  and they are orthogonal 
to the $\uuu_0 =\sqrt{\ddd}$
vector, for the coordinates of $\x_j$ the following  relations hold:
\begin{equation}\label{1}
\sum_{i=1}^n d_i x_{ji} =0,\, 
\sum_{i=1}^n d_i x_{ji}^2 =1 
\, (j=1,\dots ,k-1), \,  
\sum_{i=1}^n d_i x_{ji} x_{li} =0 \, (j\ne l).
\end{equation}
Now we will find a vector $\y =(y_1 ,\dots ,y_n )^T$ such that for it, the
conditions
\begin{equation}\label{2}
\sum_{i=1}^n d_i y_i =0 
\end{equation}
and
\begin{equation}\label{3}
\sum_{i=1}^n d_i x_{ji} y_i =0 , \quad j=1, \dots ,k-1  
\end{equation}
hold. We are looking for $\y$ in the following form:
\begin{equation}\label{4}
y_i:= \sum_{j=1}^{k-1} |x_{ji} - a_j | - b , \quad i=1,\dots ,n,
\end{equation}
where  $a_1 , \dots ,a_{k-1}$ and $b$ are appropriate real numbers.
 We will show that there exist such real numbers that
$y_i$'s defined by them satisfy conditions (\ref{2}) and (\ref{3}).

Indeed, when we already have $a_1 ,\dots ,a_{k-1}$, 
the above conditions together with $\sum_{i=1}^n d_i =1$ yield
\begin{equation}\label{5}
b=\sum_{j=1}^{k-1} b_j , \quad b_j =\sum_{i=1}^n d_i |x_{ji}-a_j | , \quad
 j=1,\dots ,k-1 .
\end{equation}
With this choice of $b$, the fulfillment of (\ref{3}) means that
for $j=1,\dots ,k-1$:
$$
\sum_{i=1}^n d_i x_{ji} y_i = \sum_{i=1}^n d_i x_{ji}
 \left( \sum_{l=1}^{k-1}  |x_{li} -a_l | -b_l \right)  = 0 .
$$
But~\eqref{1} implies that 
$$
 \sum_{i=1}^n d_i x_{ji} b_l  = 0 
$$
for $l=1,\dots ,k-1$. This provides the following system of 
equations for $a_1 ,\dots ,a_{k-1}$:
\begin{equation}\label{fontos}
 f_j =\sum_{i=1}^n d_i x_{ji} \sum_{l=1}^{k-1}  |x_{li} -a_l |  = 0 , \quad
  j=1,\dots ,k-1.
\end{equation}

We are looking for the root of the $f=(f_1 ,\dots ,f_{k-1} ): \R^{k-1} \to
\R^{k-1}$ function of stepwise linear coordinate functions.  
To prove that it has a root, we will use the multi-dimensional generalization 
of the Bolzano theorem:
a continuous map between two normed metric spaces of the same dimensions 
takes a connected set into a connected one.
Because of symmetry considerations, the range
contains the origin, see Appendix A for further details.

Now let us define the two cluster centers for the $j$th coordinates by
$$
 c_{j1} = a_j-b_j \quad \textrm{and} \quad  c_{j2} = a_j +b_j .
$$
Observe that 
$$
 |x_{ji} -a_j |-b_j  =\left\{
 \begin{array}{ll}
   c_{j1} -x_{ji} & \mbox{if} \quad x_{ji} <a_j   \\
   x_{ji} -c_{j2}   & \mbox{if} \quad x_{ji} \ge a_j  ;
 \end{array}
\right.
$$
therefore,
\begin{equation}\label{6}
|x_{ji}-a_j |-b_j =\min \{|x_{ji} -c_{j1} |,|x_{ji} -c_{j2} | \}  
\end{equation}
holds for $i=1,\dots ,n$; $j=1,\dots ,k-1$.  
For $j=1,\dots ,k-1$ they form $2^{k-1}$ centers in $k-1$
dimensions. 

Let
$$
\sigma^2(\y)=\sum_{i=1}^n d_i y_i^2
$$
be the variance of the coordinates of $\y$
with respect to the discrete measure $d_1 ,\dots ,d_n$.  
Due to (\ref{6}),   $\sigma^2(\y) \ge  { S}_{2^{k-1}}^2(\X_{k-1}^* )$.
Define the vector $\z \in \R^n$ of the following coordinates:
$$  z_i=\frac{y_i }{\sigma(\y)}, \quad i=1,\dots ,n ;
$$
obviously, $\sum_{i=1}^n d_i z_i^2 =1$.
Then
$$
\max _{i \ne m } \frac{|z_i-z_m|}{ \sum_{j=1}^{k-1} |x_{ji} -x_{jm} | }
\le \frac1{\sigma(\y)},  
$$
since due to the definition of $y_i$, the relation
$$
|y_i-y_m | \le \sum_{j=1}^{k-1} |x_{ji} -x_{jm} | \quad (i\ne m)
$$
holds. 

Let $\X_k = (\X_{k-1}^* , \z ) =(\x_1, \dots \x_{k-1} , \z )$ be 
$n\times k$ matrix, 
containing valid $k$-dimensional representatives $\rrr_1 ,\dots ,\rrr_n$ of the 
vertices in its rows;
whereas the $n\times k$ matrix
$\X_{k}^* =(\x_1 , \dots \x_{k-1} , \DD^{-1/2} \uuu_k )$
contains the optimal $k$-dimensional representatives 
in its rows. Observe that they differ only in their last coordinates.
Let $\rrr_i^*$ denote the vector comprised of the first $k-1$ coordinates of
$\rrr_i$, $i=1,\dots ,n$. These are optimal $(k-1)$-dimensional representatives
of the vertices.

By the optimality of the $k$-dimensional representation and using
Equation~\eqref{tr},   
$$
\begin{aligned}
\frac{\lambda_1+\dots +\lambda_k}{\lambda_1 +\dots \lambda_{k-1} } &=
\frac{\tr ({\X_k^*}^T \LL  \X_k^*)}{\tr ( {\X_{k-1}^*}^T \LL  
\X_{k-1}^*)} \le
\frac{\tr (\X_k^T \LL  \X_k )}{\tr({\X_{k-1}^*}^T \LL\X_{k-1}^* )}\\
& =
\frac{\sum_{i=1}^{n-1} \sum_{m=i+1}^n w_{im} \| \rrr_i-\rrr_m \|^2}
{\sum_{i=1}^{n-1} \sum_{m=i+1}^n w_{im} \| \rrr_i^* -\rrr_m^* \|^2}  \\
&=\frac{\sum_{i=1}^{n-1} \sum_{m=i+1}^n w_{im}
\left [\, \| \rrr_i^*-\rrr_m^* \|^2 + (z_i -z_m )^2 \,\right ] }
{\sum_{i=1}^{n-1} \sum_{m=i+1}^n w_{im} \| \rrr_i^* -\rrr_m^* \|^2} \\
&\le 1+\max _{ i\ne m }
\frac{(z_i-z_m)^2}{ \| \rrr_i^* -\rrr_m^* \|^2 }
\le 1 + \frac1{\sigma^2(\y)} \le 1+ \frac1{ S_{2^{k-1}}^2(\X_{k-1}^* ) },
\end{aligned}
$$
which -- by subtracting 1 from both the left- and right-hand sides and 
taking the reciprocals -- finishes the proof.
\end{proof}

Note that only if $\lambda_{k-1} <\lambda_k$, $\uuu_k$ and $\x_k$ are 
not in the subspace spanned by $\uuu_1 ,\dots ,\uuu_{k-1}$.  
Theorem~\ref{gap} indicates the following clustering property of the 
$(k-1)$th and $k$th smallest normalized Laplacian eigenvalues:
the greater the gap between them,
the better the optimal $k$-dimensional 
representatives of the vertices can be classified into $2^{k-1}$ clusters.

\begin{figure}
    \centering
    \includegraphics[scale=0.2]{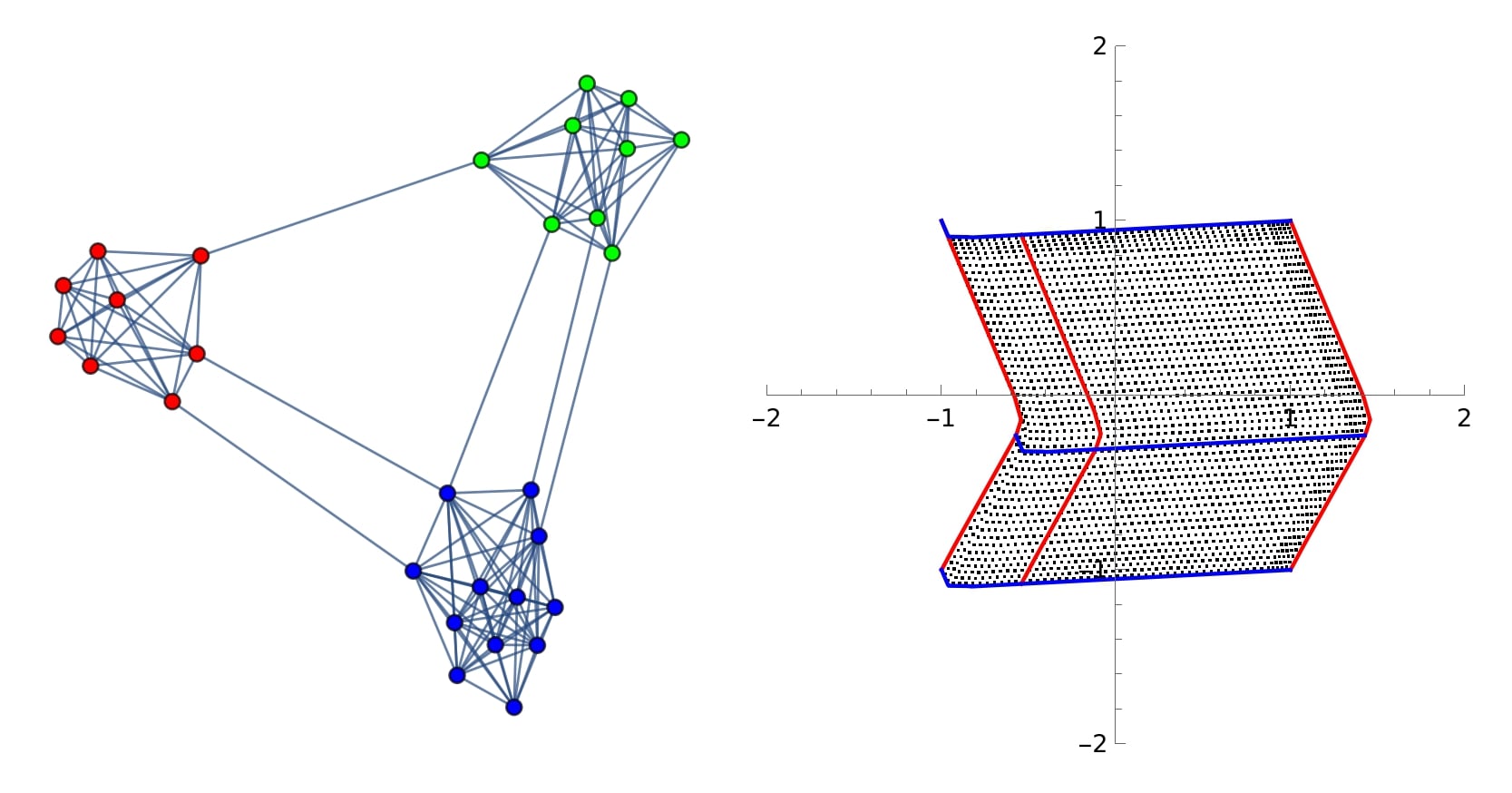}
    \caption{A graph with three well separated clusters, and the image of its Fiedler-carpet.}
    \label{fig:3clusters}
\end{figure}

Figure~\ref{fig:3clusters} shows a graph 
\footnote{This graph can be represented via the \texttt{graph6} format using the following string (see \url{http://users.cecs.anu.edu.au/~bdm/data/formats.txt} for more information):
$\backslash \backslash \sim \sim \sim \sim \{$???@$\_$F?N?n$\_$FwB$\sim $?N$\{$?ng@$\sim $w@$\sim \{$????C?G??@a??F???${}^{\wedge}$???N$\_$??FW??@$\sim $??CN$\{$}
with three well separated clusters, 
as well as the image of the corresponding map $f : \R^2 \to \R^2$ as used in 
the proof of Theorem~\ref{gap}. 

\begin{figure}
    \centering
    \includegraphics[scale=0.2]{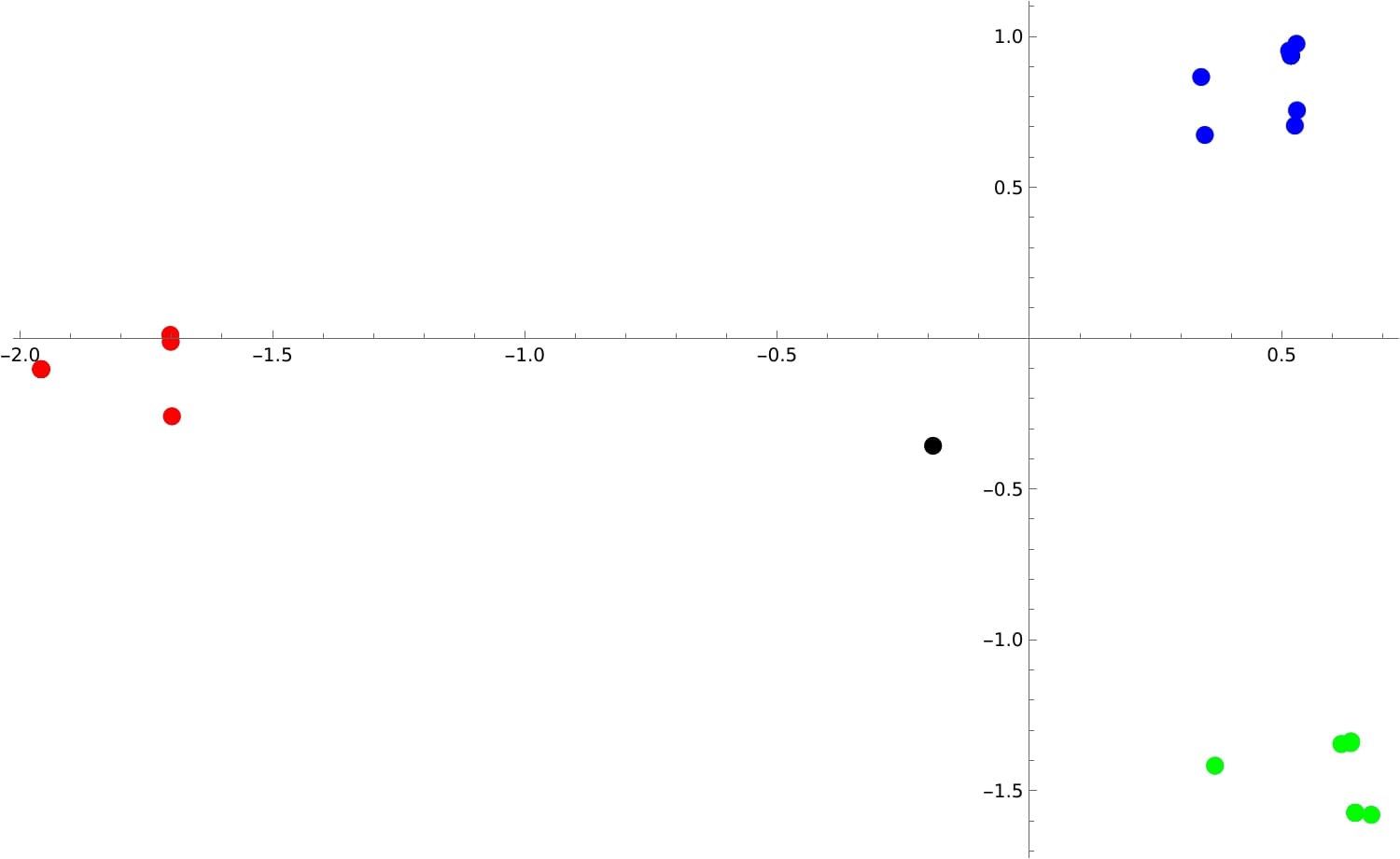}
    \caption{2-dimensional vertex-representatives of the graph of 
Figure~\ref{fig:3clusters}. 
The black point denotes the approximate root of $f$.}
    \label{fig:quadraticembedding_withroot}
\end{figure}

The image of $f$ contains the origin, and we can find by inspection a pair $(a_1,a_2)$ for which $f(a_1,a_2)$ is approximately zero; namely, choosing $a_1 = -0.19099$ and $a_2 = -0.35688$ gives $f(a_1,a_2) \approx (-0.000002, 0.00000001)$.
We compute a solution of the quadratic embedding problem in $\R^2$, where 
the 2-dimensional representative of vertex $i$ is the point $(x_{1i}, x_{2i})$
for $i=1,\dots n$, and the black point denotes the approximate root of $f$,
see Figure~\ref{fig:quadraticembedding_withroot}.

\vskip0.5cm
Some remarks are in order:
\begin{itemize}
\item
Theorem~\ref{gap} is the generalization of Theorem 2.2.3 of~\cite{Bolla3}.
There, in the $k=2$ case, $2^{k-1} =k$, but in general, the number of
clusters is much larger than that of the relevant eigenvectors.
 

\item
The statement of the theorem has relevance, since for any $k>0$, the relation
$\sum_{j=1}^{k-1}\lambda_j <(k-1) \lambda_k$ holds; but with analysis of variance
considerations, $S_{2^{k-1}}^2 \le k-1$ also holds. In particular, when $k=2$,
only one eigenvector is used for the representation. The total inertia
of the coordinates is 1, and it can be divided into the sum of nonnegative
within- and between-cluster inertias. The within-cluster inertia is the sum
of the inner variances of the two clusters, which is $S_2^2$, so it is at 
most 1. (Here the variances are calculated with respect to the discrete
distribution $d_1 ,\dots ,d_n$.) 
When $(k-1)$-dimensional representatives are used, then the total inertia is
$\tr (\X_{k-1}^T \DD \X_{k-1} ) =\tr (\I_{k-1}) =k-1$, and the sum of the inner
variances is again at most $k-1$, but it is further  bounded with
$\frac{\sum_{j=1}^{k-1}\lambda_j}{\lambda_k }$ by Theorem~\ref{gap}. 

\item
The vector $\uuu_1$ is called Fiedler-vector of the (non-normalized) Laplacian.
Here we use the first $k-1$ transformed eigenvectors of the normalized
Laplacian together, that we call \textit{Fiedler-carpet}.



\end{itemize}

\section{Applications}\label{appl}

We investigated the international migrant stock by the country of origin and destination in the years 2015 and 2019. The focus in on 41 European countries 
plus the United States of America and Canada. The data\footnote{United Nations, Department of Economic and Social Affairs, Population Division, International Migrant Stock 2019 (\url{https://www.un.org/en/development/desa/population/migration/data/estimates2/estimates19.asp)}.} is based on the official registered migrants numbers, where columns and rows are corresponding to the country of origin and the country of destination, respectively.

Here the quadratic, but not symmetric edge-weight matrix contains weights
of bidirected edges (the diagonal is zero): the $i,j$ entry is the number
of persons going $j\to i$.  Via SVD of the normalized table,
we can find emigration (column) and immigration (row) clusters, 
between which the migration is the best homogeneous
(in terms of discrepancy). 

Both for the 2015 and 2019 data, there was a gap after four 
non-trivial singular 
values, and therefore the corresponding four singular vector pairs were used to 
find five emigration and immigration trait clusters. 
\begin{itemize}
	\item Singular values, 2015:
	\par  1, 0.79098, 0.71857, 0.67213, 0.56862, 0.45293, 0.40896, 0.38178, 0.36325,  0.34785, 0.32648, 0.31769, 0.2996,  0.27927, 0.26566, 0.24718, 0.22638, 0.20632,  0.18349, 0.1651,  0.14384, 0.1359,  0.12721, 0.12092, 0.11816, 0.10374, 0.09545,  0.08278, 0.0738,  0.06371, 0.05673, 0.04553, 0.03488, 0.03107, 0.02967, 0.02693,  0.01557, 0.00788, 0.00584, 0.00519, 0.00191, 0.0017,  0.00099.

	\item Singular values, 2019:
	\par  1, 0.77844, 0.70989, 0.65059, 0.55122, 0.43612, 0.39512, 0.36194, 0.3558,  0.33882, 0.32174, 0.30719, 0.29601, 0.28181, 0.26865, 0.259,   0.22421, 0.1917,  0.17988, 0.1516,  0.13671, 0.13243, 0.12397, 0.11542, 0.10598, 0.09216, 0.08889,  0.07958, 0.06835, 0.06154, 0.05377, 0.04412, 0.03436, 0.03124, 0.02899, 0.02745,  0.01507, 0.00814, 0.00619, 0.0051,  0.00216, 0.00129, 0.00089.
\end{itemize}

\begin{table}[t]
	\begin{minipage}{0.55\textwidth}
		
		\centering
		\begin{tabular}{|c|p{3.8cm}|}
			\hline
			Cluster \#  & Emigration Countries \\
			\hline
			1  & Austria, Belgium, Bulgaria, Canada, Czechia, Denmark, Finland, France, Germany, Greece, Hungary, Iceland, Ireland, Italy, Latvia, Liechtenstein, Lithuania, Luxembourg, Malta, Monaco, Netherlands, North Macedonia, Norway, Poland, Portugal, Romania, Serbia, Slovakia, Spain, Sweden, Switzerland, United Kingdom, United States of America \\
			\hline
			2  & Bosnia and Herzegovina, Croatia, Montenegro, Slovenia \\
			\hline
			3  & Albania \\
			\hline
			4  & Belarus, Estonia, Republic of Moldova, Ukraine \\
			\hline
			5  & Russian Federation \\ \hline
		\end{tabular}
		\caption{Country memberships of emigration trait clusters, 2015.}
		\label{tab:2015_row}	
	\end{minipage}
	\hfil
	\begin{minipage}{0.55\textwidth}
		\centering
		\begin{tabular}{|c|p{4cm}|}
			\hline
			Cluster \#  & Immigration Countries \\
			\hline
			1  & Albania, Austria, Belgium, Bulgaria, Canada, Czechia, Denmark, Finland, France, Germany, Hungary, Iceland, Ireland, Liechtenstein, Luxembourg, Malta, Monaco, Netherlands, Norway, Portugal, Romania, Slovakia, Spain, Sweden, Switzerland, United Kingdom, United States of America 
			\\
			\hline
			2  & Greece, Italy, North Macedonia \\
			\hline
			3  & Bosnia and Herzegovina, Croatia, Montenegro, Serbia, Slovenia \\
			\hline
			4  & Belarus, Estonia, Latvia, Lithuania, Ukraine \\
			\hline
			5  & Poland, Republic of Moldova, Russian Federation
			\\ 		\hline
		\end{tabular}		
		\caption{Country memberships of immigration trait clusters, 2015.}
		\label{tab:2015_col}
	\end{minipage}
\end{table}

\begin{table}[t]
	\begin{minipage}{0.55\textwidth}
		\centering
\begin{tabular}{|c|p{4cm}|}
	\hline
	Cluster \# & Emigration Countries \\
	\hline
	1  & Austria, Belgium, Bulgaria, Canada, Czechia, Denmark, Finland, France, Germany, Greece, Hungary, Iceland, Ireland, Italy, Latvia, Liechtenstein, Lithuania, Luxembourg, Malta, Monaco, Netherlands, North Macedonia, Norway, Poland, Portugal, Romania, Serbia, Slovakia, Spain, Sweden, Switzerland, United Kingdom, United States of America \\
	\hline
	2  & Bosnia and Herzegovina, Croatia, Montenegro, Slovenia \\
	\hline
	3  & Albania \\
	\hline
	4  & Belarus, Estonia, Republic of Moldova, Ukraine \\
	\hline
	5  & Russian Federation \\
	\hline
\end{tabular}
\caption{Country memberships of emigration trait clusters, 2019.}
\label{tab:2019_col}
	\end{minipage}
	\hfil
	\begin{minipage}{0.55\textwidth}
	\centering
\begin{tabular}{|c|p{4.1cm}|}
	\hline
	Cluster \#  & Immigration Countries \\
	\hline
	1  & Albania, Austria, Belgium, Bulgaria, Canada, Czechia, Denmark, Finland, France, Germany, Hungary, Iceland, Ireland, Liechtenstein, Luxembourg, Malta, Monaco, Netherlands, Norway, Portugal, Romania, Slovakia, Spain, Sweden, Switzerland, United Kingdom, United States of America \\
	\hline
	2  & Bosnia and Herzegovina, Croatia, Montenegro, Serbia, Slovenia \\
	\hline
	3  & Greece, Italy, North Macedonia \\
	\hline
	4  & Poland, Republic of Moldova, Russian Federation \\
	\hline
	5  & Belarus, Estonia, Latvia, Lithuania, Ukraine \\
	\hline
\end{tabular}
\caption{Country memberships of immigration trait clusters, 2019.}
\label{tab:2019_row}
	\end{minipage}
\end{table}



\begin{figure}[htp]
	\begin{minipage}{.5\textwidth}
		\centering
		\includegraphics[width=6.5cm]{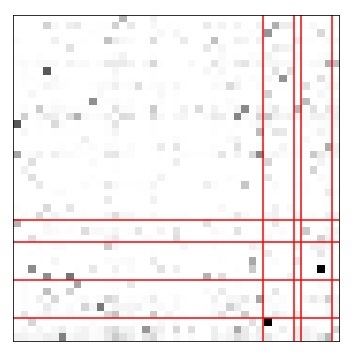}
		\subcaption{2015}
		\label{fig:2015_cross}
	\end{minipage}		
	\hfil
	\begin{minipage}{.5\textwidth} 			
		\centering
		\includegraphics[width=6.5cm]{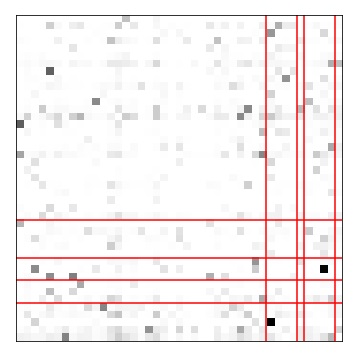}
		\subcaption{2019}
		\label{fig:2019_cross}
	\end{minipage}			
	\caption{Immigration--emigration cluster-pairs with the countries 
rearranged by their cluster memberships. The frequency counts are represented 
with light to dark squares.}
\end{figure}	

Emigration and immigration trait clusters for 2015 are shown in
Table~\ref{tab:2015_col} and Table~\ref{tab:2015_row}, whereas those for
2019, in Table~\ref{tab:2019_col} and Table~\ref{tab:2019_row}.
Figures~\ref{fig:2015_cross} and~\ref{fig:2019_cross} illustrate the
immigration--emigration cluster-pairs with the countries rearranged by
their cluster memberships. The frequency counts are represented with 
light to dark squares.

In 2015, by $\chi^2$ test we found the smallest discrepancy between the 
emigration trait cluster number 2  and immigration trait cluster number 4,
i.e. the sub-table formed by them was close to an independent table of rank 1.
The clusters were similar in the two years; 
in 2019, the smallest discrepancy was found between the 
emigration trait cluster number 2  and immigration trait cluster number 5,
i.e. between the Balcanian and Baltic countries in both years.

Correspondence analysis results with cluster memberships are found in 
Appendix B.

\section{Summary}\label{sum}

Spectral clustering is for clustering data points or the vertices of a graph
based on combinatorial criteria with spectral relaxation.
Here we  generalize spectral clustering in several ways:
\begin{itemize}
\item
Instead of simple or edge-weighed graphs consider directed graphs and
rectangular arrays of
  nonnegative entries. For these, the method of correspondence analysis
  gives low-dimensional representation of the row and column items by
  means of SVD of the normalized table. This is applied to real-life data.

\item
Instead of the usual multiway minimum cut objective, we consider
  discrepancy minimization, for which a near optimal solution is given
  by the $k$-means algorithm applied to the row and column representatives,
  where the number of clusters $(k)$ is concluded from the gap in the singular
  spectrum. There are theoretical results supporting this.
  Eventually, the near optimal clusters can be refined to decrease discrepancy.
\item The number of clusters can be larger than the number of eigenvectors
entered in the classification. The two are the same only in case of
generalized random or quasirandom graphs, see~\cite{Bolla20}.
In Theorem~\ref{gap}, 
an exact estimate for the sum of the inner variances of $2^{k-1}$ clusters is
given by means of the spectral gap between the $(k-1)$th and $k$th
smallest normalized Laplacian eigenvalues by using the Fiedler-carpet
of the corresponding eigenvectors.
\end{itemize}

\section*{Acknowledgements}

The research was done under the auspices of the Budapest Semesters in 
Mathematics program, in the framework of an undergraduate online research 
course in the fall semester
2021, with the participation of US undergraduate students.
Also, Fatma Abdelkhalek's work is funded by a scholarship under the 
Stipendium Hungaricum program between Egypt and Hungary.
The paper is dedicated to G\'abor Tusn\'ady for his 80th birthday,
with whom the first author of this paper posed a conjecture in~\cite{Bolla3}.
Here the conjecture is proved in the $k>2$ case, albeit with much more
clusters than the number of the eigenvalues preceding the spectral gap.

\newpage

\section*{Appendix A}\label{appA}

In the $f:\R \to  \R$ case: $f=f_1$, $a=a_1$ and with the notation
$A =\min_i x_{1i}$, $B=\max_i x_{1i}$ for
$$
 f(a)= \sum_{i=1}^n d_i x_{1i} |x_{1i} -a |
$$
we have that $f(A)=1$ and $f(B)=-1$. As $f$ is continuous, it must have a root
in $(A,B)$, by the Bolzano theorem. 
Also note that this root of $f$ is
around the median of the coordinates of $\x_1$ with respect to the discrete
measure $d_1 ,\dots ,d_n$.

Then consider the $f= (f_1 ,f_2 ) :\R^2 \to  \R^2$ case. 
The coordinate functions are stepwise linear,
continuous functions. The system of equations is
\begin{equation}\label{f}
\begin{aligned}
 f_1 (a_1 ,a_2 ) &= \sum_{i=1}^n d_i x_{1i} |x_{1i} -a_1 |  +
                 \sum_{i=1}^n d_i x_{1i} |x_{2i} -a_2 |= 0 , \\
 f_2 (a_1 ,a_2 ) &= \sum_{i=1}^n d_i x_{2i} |x_{1i} -a_1 |  +
                 \sum_{i=1}^n d_i x_{2i} |x_{2i} -a_2 |= 0 .
\end{aligned}
\end{equation}
 With the notation $A =\min_i x_{1i}$, $B=\max_i x_{1i}$,
$C =\min_i x_{2i}$, $D=\max_i x_{2i}$, where $A<0,B>0,C<0,D>0$,
$$
f(A, C) =(1,1), \, f(B, D) =(-1,-1), \, f(A, D) =(1,-1), \,
 f(B, C) =(-1,1) .
$$
Furthermore,
$f(x, y) =(1,1)$ if $x\le A, \, y\le C$; $f(x, y) =(-1,-1)$ if $x\ge B, \,
y\ge D$; $f(x, y) =(1,-1)$ if $x\le A, \, y\ge D$;  and
$f(B, C) =(-1,1)$ if $x\ge B ,\, y\le C$.

We want to show that $f$ has a root within the rectangle
$[A ,B]\times [C,D]$. By the multivariate version of the Bolzano theorem,
the $f$-map of this rectangle is a connected region in $\R^2$ that contains the
points $(1,1), (-1,-1), (1,-1), (-1,1)$ as `corners'.
We will show that it contains $(0,0)$ too. 

Note that together with $\uuu_j$, the vector $-\uuu_j$ is also a unit-norm
eigenvector, so instead of $\x_j$ we can as well use $-\x_j$ for $j=1,2$, 
that gives 4 possible domains of $f$: next to the rectangle
$[A ,B]\times [C,D]$, the rectangles $[-B ,-A]\times [C,D]$, 
$[A ,B]\times [-D,-C]$, and
$[-B ,- A]\times [-D,-C]$ are also closed, bounded regions,
and the $f$-images of them show symmetry with respect to the coordinate axes.
Therefore, it suffices 
to prove that the map of the union of them contains the origin.
In Section~\ref{pre}, we saw that neither the objective
function ($Q_k$), nor the clustering is affected by the orientation of the 
eigenvectors, so the orientation is not denoted in the sequel.
Also notice that with counter-orienting $\uuu_1$ or/and $\uuu_2$: if
$(a_1, a_2)$ is a root of $f$, then $-a_1$ instead of $a_1$ or/and 
$-a_2$ instead of $a_2$ will result in a root of $f$ too.

The images are closed, bounded  regions (usually not rectangles),
but we will show that the opposite sides of  them are parallel curves
and sandwich the $f_1$ and $f_2$ axes, respectively.
As the $f$-values sweep the region between these boundaries, 
the total range should contain the origin. 
At the end, we take the orientations into consideration
to arrive to the expected result.  

Now the above, below, right, and left boundaries are investigated.

\begin{itemize}
\item
\textbf{Above:} Consider the boundary curve between (-1,1) and (1,1). 
Along that, $a_2 =C$ and $A<a_1 <B$. 
Let $H:= \{ i: \, x_{1i} >a_1 \}$. 
Then $H\neq \emptyset$ and ${\bar H} \neq \emptyset$; further,
\begin{equation}\label{contur}
\begin{aligned}
f_2 (a_1 ,C) &=\sum_{i=1}^n d_i x_{2i} |x_{1i} -a_1 | 
             +\sum_{i=1}^n d_i x_{2i} (x_{2i} -C ) = \\
&=\sum_{i\in H} d_i x_{2i} (x_{1i} -a_1 ) -
 \sum_{i\in {\bar H}} d_i x_{2i} (x_{1i} -a_1 ) + 1 \\
&=2\sum_{i\in H} d_i x_{2i} (x_{1i} -a_1 ) +1 \\
&=2\sum_{i\in {\bar H}} d_i x_{2i} (a_1 - x_{1i} ) +1 ,
\end{aligned}
\end{equation}
where we intensively used conditions~\eqref{1}.

\item
\textbf{Below:} Consider the boundary curve between (-1,-1) and (1,-1). 
Along that, $a_2 =D$ and $A<a_1 <B$. 
Then
\begin{equation}\label{contur1}
\begin{aligned}
f_2 (a_1 ,D) &=\sum_{i=1}^n d_i x_{2i} |x_{1i} -a_1 | 
             -\sum_{i=1}^n d_i x_{2i} (x_{2i} -D ) = \\
&=\sum_{i\in H} d_i x_{2i} (x_{1i} -a_1 ) -
 \sum_{i\in {\bar H}} d_i x_{2i} (x_{1i} -a_1 ) -1 = \\
&=2\sum_{i\in H} d_i x_{2i} (x_{1i} -a_1 ) -1 \\
&=2\sum_{i\in {\bar H}} d_i x_{2i} (a_1 - x_{1i} ) -1 .
\end{aligned}
\end{equation}

\item
\textbf{Between (horizontally):} Consider the case when $a_2 =u\in (C,D)$ fixed 
and $A<a_1 <B$.  Then
$$
\begin{aligned}
f_2 (a_1 ,u ) &=\sum_{i=1}^n d_i x_{2i} |x_{1i} -a_1 | 
             +\sum_{i=1}^n d_i x_{2i} | x_{2i} -u | = \\
&=\sum_{i\in H} d_i x_{2i} (x_{1i} -a_1 ) -
 \sum_{i\in {\bar H}} d_i x_{2i} (x_{1i} -a_1 ) + \\ 
&+ \sum_{i:\, x_{2i} >u } d_i x_{2i} (x_{2i} -u )+
 \sum_{i:\, x_{2i}\le u } d_i x_{2i} (u-x_{2i}  ) = \\
&=2\sum_{i\in H} d_i x_{2i} (x_{1i} -a_1 ) -1+ 
 2 \sum_{i:\, x_{2i} >u } d_i x_{2i}^2 -2u \sum_{i:\, x_{2i} >u } d_i x_{2i} .
\end{aligned}
$$

\end{itemize}
So the $f_2 (a_1 ,u)$ arcs are all  parallel 
to the boundary curves $f_2 (a_1 ,C)$ and $f_2 (a_1 ,D)$ and to each other.
In particular,
$$
f_2 (a_1 ,0 )  =
 2\sum_{i\in H} d_i x_{2i} (x_{1i} -a_1 ) -1+ 
 2 \sum_{i:\, x_{2i} >0 } d_i x_{2i}^2 .
$$
This arc is either closer to the Above or the Below curve (which are in
distance 2 from each other), depending on
whether $\sum_{i: \, x_{2i} >0} d_i x_{2i}^2$ is less or greater than 
$\frac12$, but is strictly positive by condition~\eqref{1}. 
For this reason, if $u$ and $u'$ are `close' to 0, then
the $f_2 (a_1 ,u )$ and $f_2 (a_1 ,u' )$ arcs are not identical, 
otherwise it can happen that for some $u\ne u'$:
\begin{equation}\label{bad}
 \sum_{i:\, x_{2i} >u } d_i x_{2i}^2 -u \sum_{i:\, x_{2i} >u } d_i x_{2i} =
 \sum_{i:\, x_{2i} >u' } d_i x_{2i}^2 -u' \sum_{i:\, x_{2i} >u' } d_i x_{2i} .
\end{equation}

\vskip0.2cm
The same is true vertically. 
\begin{itemize}
\item
\textbf{Right:} 
Consider the boundary curve between (1,-1) and (1,1). 
Along that, $a_1 =A$ and $C<a_2 <D$. 
Let $F:= \{ i: \, x_{2i} >a_2 \}$. 
Then $F\neq \emptyset$ and ${\bar F} \ne \emptyset$; further,
\begin{equation}\label{right}
\begin{aligned}
f_1 (A ,a_2 ) &=\sum_{i=1}^n d_i x_{1i} (x_{1i} -A ) 
             -\sum_{i=1}^n d_i x_{1i} |x_{2i} -a_2 | = \\
&=1+ \sum_{i\in F} d_i x_{1i} (x_{2i} -a_2 ) -
 \sum_{i\in {\bar F}} d_i x_{1i} (x_{2i} -a_2 ) = \\
&=1+2\sum_{i\in F} d_i x_{1i} (x_{2i} -a_2 ) = \\
&=1+2\sum_{i\in {\bar F}} d_i x_{1i} (a_2 - x_{2i} ) .
\end{aligned}
\end{equation}

\item
\textbf{Left:} As from the left,
consider the boundary curve between (-1,-1) and (-1,1). 
Along that, $a_1 =B$ and $C<a_2 <D$. Then
\begin{equation}\label{left}
\begin{aligned}
f_1 (B ,a_2 ) 
&=-1+2\sum_{i\in F} d_i x_{1i} (x_{2i} -a_2 ) = \\
&=-1+2\sum_{i\in {\bar F}} d_i x_{1i} (a_2 - x_{2i} )  .
\end{aligned} 
\end{equation}

\item
\textbf{Between (vertically):} Consider the case when $a_1 =v \in (A,B)$ fixed 
and $C<a_2 <D$.  Then
$$
f_1 (v ,a_2 ) = 
 1+2\sum_{i\in F} d_i x_{1i} (x_{2i} -a_2 )  -1+ 
 2 \sum_{i:\, x_{1i} >v } d_i x_{1i}^2 -2v \sum_{i:\, x_{1i} >v } d_i x_{1i} .
$$
\end{itemize}
So the $f_1 (v, a_2)$ arcs are all  parallel 
to the boundary curves $f_1 (A ,a_2)$ and $f_1 (B, a_2)$ and to each other.
In particular,
$$
f_1 (0 ,a_2 )  =
 2\sum_{i\in F} d_i x_{1i} (x_{2i} -a_2 ) -1+ 
 2 \sum_{i:\, x_{1i} >0 } d_i x_{1i}^2 .
$$
This arc is either closer to the Right or the Left curve (which are in
distance 2 from each other), depending on
whether $\sum_{i: \, x_{1i} >0} d_i x_{1i}^2$ is less or greater than 
$\frac12$,  but is strictly positive by condition~\eqref{1}. 
For this reason, if $v$ and $v'$ are `close' to 0, then
the $f_1 (v, a_2 )$ and $f_1 (v', a_2 )$ arcs are not identical, 
otherwise it can happen that for some $v\ne v'$:
\begin{equation}\label{bad1}
 \sum_{i:\, x_{1i} >v } d_i x_{1i}^2 -v \sum_{i:\, x_{1i} >v } d_i x_{1i} =
 \sum_{i:\, x_{`i} >v' } d_i x_{1i}^2 -v' \sum_{i:\, x_{1i} >v' } d_i x_{1i} .
\end{equation}

Therefore, any grid on the rectangle of the domain (its horizontal and vertical
lines parallel to the $a_1$ and $a_2$ axes) is  mapped by $f$ onto a 
lattice 
with horizontal and vertical,  parallel arcs. This proves that $f$
is one-to-one whenever these arcs are not identical. 
The possible inconvenient phenomenon, when both Equations~\eqref{bad}
and~\eqref{bad1} hold for some $u\ne u'$ and $v\ne v'$ pairs, is 
experienced at the dark parts of Figure~\ref{fig:fourientations} 
near the boundaries. However, $f$ is injective in the
neighborhood of the origin that does not contain any eigenvector
coordinates (because there it has purely linear coordinate functions). 
This also depends on the underlying graph: if it shows
symmetries, then its weighted Laplacian has  multiple eigenvalues and/or
multiple coordinates of the eigenvectors that may cause complications.


To prove that the Above--Below boundaries sandwich the $f_1$ axis and
the Right--Left boundaries sandwich the $f_2$ axis, respectively, the following
investigations are made. Because the investigations are of similar vein,  
only the first  of them will be discussed  in details. We distinguish
between eight cases (denoted by an acronym), depending on, 
which half of which boundary is considered.
The estimates are supported by specific orientations of the unit norm 
eigenvectors. If $\uuu_1$ is oriented so that for the coordinates of
of $\x_1 =\DD^{-1/2} \uuu_1$ 
$$
 \sum_{i: \, x_{1i} >0} d_i x_{1i}^2 <\frac12 
$$
holds, then it is called \textit{positive orientation}, whereas, 
the opposite is negative.
 Likewise, he orientation of  $\uuu_2$ is positive if
 for the coordinates of of $\x_2 =\DD^{-1/2} \uuu_2$ 
\begin{equation}\label{kenyes}
 \sum_{i: \, x_{2i} >0} d_i x_{2i}^2 <\frac12 
\end{equation}
holds, otherwise it is negative.

\begin{itemize}
\item[{\textbf{AL}}] (Above boundary, Left half), when $a_1 >0$:
using the last but one line of  Equation~\eqref{contur},
$$
\begin{aligned}
f_2 (a_1 ,C) &=2\sum_{i\in H} d_i x_{2i} (x_{1i} -a_1 ) +1  \\
&=2\sum_{i\in H, \, x_{2i} >0} d_i x_{2i} (x_{1i} -a_1 ) +
  2\sum_{i\in H, \, x_{2i} \le 0} d_i x_{2i} (x_{1i} -a_1 ) + 1  \\
&\ge 2\sum_{i\in H, \, x_{2i} \le 0} d_i x_{2i} (x_{1i} -a_1 ) + 1 .
\end{aligned}
$$
To prove that  $f_2 (a_1 ,C) \ge 0$ it suffices to prove that
$$
  \sum_{i\in H, \, x_{2i} \le 0} d_i x_{2i} (x_{1i} -a_1 ) \ge -\frac12 ,
$$
which is equivalent to
$$
  \sum_{i\in H, \, x_{2i} \le 0} d_i (-x_{2i} ) (x_{1i} -a_1 ) \le \frac12 .
$$
We  use the Cauchy--Schwarz inequality by keeping in mind that 
$x_{1i} -a_1 >0$
$(i\in H)$ and because of $a_1 >0$, $x_{1i} -a_1 <x_{1i}$. Therefore, 
$$
 [\sum_{i\in H \, x_{2i} \le 0} (\sqrt{d_i} (-x_{2i})) (\sqrt{d_i} (x_{1i} -a_1 ) )]^2 
\le [ \sum_{i\in H \, x_{2i} \le 0} d_i x_{2i}^2 ] 
    [ \sum_{i\in H \, x_{2i} \le 0} d_i x_{1i}^2  ] \le 
 \frac12 \frac12
$$
holds true if $\uuu_1$ is positively and $\uuu_2$ is negatively oriented.

\item[{\textbf{AR}}] (Above boundary, Right half), when $a_1 \le 0$: 
using the last line of Equation~\eqref{contur},
to prove that  $f_2 (a_1 ,C) \ge 0$ it suffices to prove that
$$
  \sum_{i\in {\bar H}, \, x_{2i} \le 0} d_i x_{2i} (a_1 -x_{1i} ) \ge -\frac12 ,
$$
which is  equivalent to
$$
  \sum_{i\in {\bar H}, \, x_{2i} \le 0} d_i (-x_{2i} ) (a_1 - x_{1i} ) \le \frac12 .
$$
We  again use the Cauchy--Schwarz inequality by keeping in mind that 
$a_1 -x_{1i} \ge 0$
$(i\in {\bar H} )$ and $a_1 -x_{1i} = -x_{i1} - (-a_1 ) \le -x_{1i}$ as now
$-a_1 \ge 0$ and $-x_{1i} > -a_1$. Therefore,
$$
 [\sum_{i\in {\bar H }, x_{2i} \le 0 } (\sqrt{d_i} (-x_{2i})) (\sqrt{d_i} 
      ( a_1 -x_{1i} ) )]^2
 \le  [ \sum_{i\in {\bar H} ,  x_{2i} \le 0 } d_i x_{2i}^2 ] 
 [ \sum_{i\in {\bar H} ,  x_{2i} \le 0 } d_i (-x_{1i})^2  ] \le 
 \frac12 \frac12
$$
holds true with negatively orienting $\uuu_1$ and negatively $\uuu_2$.

\item[{\textbf{BL}}] (Below boundary, Left half), when $a_1 > 0$: 
using the last but one line of  Equation~\eqref{contur1},
$$
\begin{aligned}
f_2 (a_1 ,D) &=2\sum_{i\in H} d_i x_{2i} (x_{1i} -a_1 ) -1  \\
&=2\sum_{i\in H, \, x_{2i} < 0} d_i x_{2i} (x_{1i} -a_1 ) +
  2\sum_{i\in H, \, x_{2i} \ge 0} d_i x_{2i} (x_{1i} -a_1 ) - 1  \\
&\le 2\sum_{i\in H, \, x_{2i} \ge 0} d_i x_{2i} (x_{1i} -a_1 ) - 1 .
\end{aligned}
$$
To prove that  $f_2 (a_1 ,D) \le 0$ it suffices to prove that
$$
  \sum_{i\in H, \, x_{2i} \ge 0} d_i x_{2i} (x_{1i} -a_1 ) \le \frac12 .
$$
We  use the Cauchy--Schwarz inequality by keeping in mind that 
$x_{1i} -a_1 >0$
$(i\in H)$ and because of $a_1 >0$, $x_{1i} -a_1 <x_{1i}$. Therefore, 
$$
 [\sum_{i\in H \, x_{2i} \ge 0} (\sqrt{d_i} x_{2i} ) (\sqrt{d_i} (x_{1i} -a_1 ) )]^2 
\le [ \sum_{i\in H \, x_{2i} \ge 0} d_i x_{2i}^2 ] 
    [ \sum_{i\in H \, x_{2i} \ge 0} d_i x_{1i}^2  ] \le 
 \frac12 \frac12
$$
holds true if $\uuu_1$ is positively and $\uuu_2$ is positively oriented.

\item[{\textbf{BR}}] (Below boundary, Right half), when $a_1 \le 0$:  
using the last line of Equation~\eqref{contur1},
to prove that  $f_2 (a_1 ,D) \le 0$ it suffices to prove that
$$
  \sum_{i\in {\bar H}, \, x_{2i} \ge 0} d_i x_{2i} (a_1 - x_{1i} ) \le \frac12 .
$$
We  again use the Cauchy--Schwarz inequality by keeping in mind that 
$a_1 -x_{1i} \ge 0$
$(i\in {\bar H} )$ and $a_1 -x_{1i} = -x_{i1} - (-a_1 ) \le -x_{1i}$ as now
$-a_1 \ge 0$ and $-x_{1i} > -a_1$. Therefore,
$$
 [\sum_{i\in {\bar H }, x_{2i} \ge 0 } (\sqrt{d_i} x_{2i}) (\sqrt{d_i} 
      ( a_1 -x_{1i} ) )]^2
 \le  [ \sum_{i\in {\bar H} ,  x_{2i} \ge 0 } d_i x_{2i}^2 ] 
 [ \sum_{i\in {\bar H} ,  x_{2i} \ge 0 } d_i (-x_{1i})^2  ] \le 
 \frac12 \frac12
$$
holds true with negatively orienting $\uuu_1$ and positively $\uuu_2$.

\item[{\textbf{RB}}] (Right boundary, Below half), when $a_2 > 0$: 
using the last but one line of  Equation~\eqref{right},
$$
f_1 (A ,a_2) 
\ge 2\sum_{i\in F, \, x_{1i} \le 0} d_i x_{1i} (x_{2i} -a_2 ) + 1 .
$$
To prove that  $f_1 (A, a_2 ) \ge 0$ it suffices to prove that
$$
  \sum_{i\in F, \, x_{1i} \le 0} d_i x_{1i} (x_{2i} -a_2 ) \ge -\frac12 ,
$$
which is equivalent to
$$
  \sum_{i\in F, \, x_{1i} \le 0} d_i (-x_{1i} ) (x_{2i} -a_2 ) \le \frac12 .
$$
By the Cauchy--Schwarz inequality,
$$
 [\sum_{i\in F \, x_{1i} \le 0} (\sqrt{d_i} (-x_{1i})) 
 (\sqrt{d_i} (x_{2i} -a_2 ) )]^2 
\le [ \sum_{i\in F \, x_{1i} \le 0} d_i x_{1i}^2 ] 
    [ \sum_{i\in F \, x_{1i} \le 0} d_i x_{2i}^2  ] \le 
 \frac12 \frac12
$$
holds true if $\uuu_1$ is negatively and $\uuu_2$ is positively oriented.

\item[{\textbf{RA}}] (Right boundary, Above half), when $a_2 \le 0$: 
using the last line of Equation~\eqref{right},
to prove that  $f_1 (A, a_2 ) \ge 0$ it suffices to prove that
$$
  \sum_{i\in {\bar F}, \, x_{1i} \le 0} d_i x_{1i} (a_2 -x_{2i} ) \ge -\frac12 ,
$$
which is  equivalent to
$$
  \sum_{i\in {\bar F}, \, x_{1i} \le 0} d_i (-x_{1i} ) (a_2 - x_{2i} ) \le \frac12 .
$$
By the Cauchy--Schwarz inequality,
$$
 [\sum_{i\in {\bar F }, x_{1i} \le 0 } (\sqrt{d_i} (-x_{1i})) (\sqrt{d_i} 
      ( a_2 -x_{2i} ) )]^2
 \le  [ \sum_{i\in {\bar F} ,  x_{1i} \le 0 } d_i x_{1i}^2 ] 
 [ \sum_{i\in {\bar F} ,  x_{1i} \le 0 } d_i (-x_{2i})^2  ] \le 
 \frac12 \frac12
$$
holds true with negatively orienting $\uuu_1$ and negatively $\uuu_2$.

\item[{\textbf{LB}}] (Left boundary, Below half), when $a_2 > 0$: 
using the last but one line of  Equation~\eqref{left},
$$
f_1 (B, a_2 ) =2\sum_{i\in F} d_i x_{1i} (x_{2i} -a_2) -1  
 \le 2\sum_{i\in F, \, x_{1i} \ge 0} d_i x_{1i} (x_{2i} -a_2 ) - 1 .
$$
To prove that  $f_1 (B, a_2 ) \le 0$ it suffices to prove that
$$
  \sum_{i\in F, \, x_{1i} \ge 0} d_i x_{1i} (x_{2i} -a_2 ) \le \frac12 ,
$$
By the Cauchy--Schwarz inequality,
$$
 [\sum_{i\in F \, x_{1i} \ge 0} (\sqrt{d_i} x_{1i}) (\sqrt{d_i} (x_{2i} -a_2 ) )]^2 
\le [ \sum_{i\in F \, x_{1i} \ge 0} d_i x_{1i}^2 ] 
    [ \sum_{i\in F \, x_{1i} \ge 0} d_i x_{2i}^2  ] \le 
 \frac12 \frac12
$$
holds true if $\uuu_1$ is positively and $\uuu_2$ is positively oriented.

\item[{\textbf{LA}}] (Left boundary, Above half), when $a_2 \le  0$:  
using the last line of Equation~\eqref{left},
to prove that  $f_1 (B, a_2 ) \le 0$ it suffices to prove that
$$
  \sum_{i\in {\bar F}, \, x_{1i} \ge 0} d_i x_{1i} (a_2 - x_{2i} ) \le \frac12 .
$$
By the Cauchy--Schwarz inequality, 
$$
 [\sum_{i\in {\bar F }, x_{1i} \ge 0 } (\sqrt{d_i} x_{1i}) (\sqrt{d_i} 
      ( a_2 -x_{2i} ) )]^2
 \le  [ \sum_{i\in {\bar F} ,  x_{1i} \ge 0 } d_i x_{1i}^2 ] 
 [ \sum_{i\in {\bar F} ,  x_{1i} \ge 0 } d_i (-x_{2i})^2  ] \le 
 \frac12 \frac12
$$
holds true with positively orienting $\uuu_1$ and negatively $\uuu_2$.
\end{itemize}

So the convenient orientation of the \textbf{AL} scenario matches that of the 
\textbf{LA} one. Similarly, the \textbf{AR-RA}, \textbf{BL-LB}, and 
\textbf{BR-RB} scenatios can be ralized with the same orientation of $\uuu_1 ,
\uuu_2$.
Namely,
the first quadrant of the domain $(a_1 \ge 0 , a_2 \ge 0)$ is mapped to the
third quadrant of the range, in which part the $(+,+)$ orientation;
the second quadrant of the domain $(a_1 \le 0 , a_2 \ge 0)$ is mapped to the
fourth quadrant of the range, in which part the $(-,+)$ orientation;
the third quadrant of the domain $(a_1 \le 0 , a_2 \le 0)$ is mapped to the
first quadrant of the range, in which part the $(-,-)$ orientation;
the fourth quadrant of the domain $(a_1 \ge 0 , a_2 \le 0)$ is mapped to the
second quadrant of the range, in which part the $(+,-)$ orientation of
$(\x_1 ,\x_2 )$ works.

Therefore, the ranges corresponding to the four quadrants intesect, and
so, the ranges under the four different orientations are connected regions
(by the multivariate analogue of the Bolzano theorem) and all contain the
`corners'  $(1,1), (-1,-1), (1,-1), (-1,1)$. Therefore, the union is also a
connected region in $\R^2$. As it is bounded from above, from below, from
the right, and from the left with curves that enclose the origin, it
should contain the origin too. Moreover, there should be an orientation and a
quadrant that contains the origin. It is not important, but note that with
any orientation, a quadrant contains the origin, which is the map of the
root we are looking for.

Alternatively, notice that in the \textbf{AL} case,
$f_2 (a_1 ,C)$ is greater than or less than 1, depending on whether the
absolute value of
$\sum_{i\in H, \, x_{2i} >0} d_i x_{2i} (x_{1i} -a_1 )$ or that of
$\sum_{i\in H, \, x_{2i} \le 0} d_i x_{2i} (x_{1i} -a_1 )$  is larger.
By the Cauchy--Schwarz inequality this happens with the positive orientation of
$\uuu_1$ and either with the positive or negative orientation of $\uuu_2$.
Simultaneously, in the \textbf{BL} case, 
$f_2 (a_1 ,D)$ greater than or less than -1, depending on whether the
absolute value of
$\sum_{i\in H, \, x_{2i} >0} d_i x_{2i} (x_{1i} -a_1 )$ or that of
$\sum_{i\in H, \, x_{2i} \le 0} d_i x_{2i} (x_{1i} -a_1 )$  is larger.
So with the positive orientation of
$\uuu_1$ and either with the positive or negative orientation of $\uuu_2$,
$f_2 (a_1 ,C )$ and $f_2 (a_1 ,D )$  are between -2 and 2, they are 
parallel, in distance 2 from each other and sandwich the $f_1$ axis.

Likewise,  the \textbf{AR,BR} estimates imply that
with the negative orientation of
$\uuu_1$ and either with the positive or negative orientation of $\uuu_2$,
$f_2 (a_1 ,C )$ and $f_2 (a_1 ,D )$  are between -2 and 2, they are 
parallel, in distance 2 from each other and sandwich the $f_1$ axis.
However, changing the orientation of $\uuu_1$ just means that 
$-a_1$ instead of $a_1$ can be in the root of $f$.

Vertically,
the \textbf{RB,LB} estimates imply that
with the positive orientation of
$\uuu_2$ and either with the positive or negative orientation of $\uuu_1$,
$f_2 (A, a_2 )$ and $f_2 (B, a_2 )$  are between -2 and 2, they are 
parallel, in distance 2 from each other and sandwich the $f_2$ axis.
Likewise,  the \textbf{AR,LA} estimates imply that
with the negative orientation of
$\uuu_2$ and either with the positive or negative orientation of $\uuu_1$,
$f_2 (A, a_2 )$ and $f_2 (B, a_2 )$  are between -2 and 2, they are 
parallel, in distance 2 from each other and sandwich the $f_2$ axis.
However, changing the orientation of $\uuu_2$ just means that 
$-a_2$ instead of $a_2$ can be in the root of $f$.

As two of these orientations match each other,   
the $f_1$ and $f_2$ axes are sandwiched, 
and so, the origin should be included.

\begin{figure}
    \centering
    \includegraphics[scale=0.2]{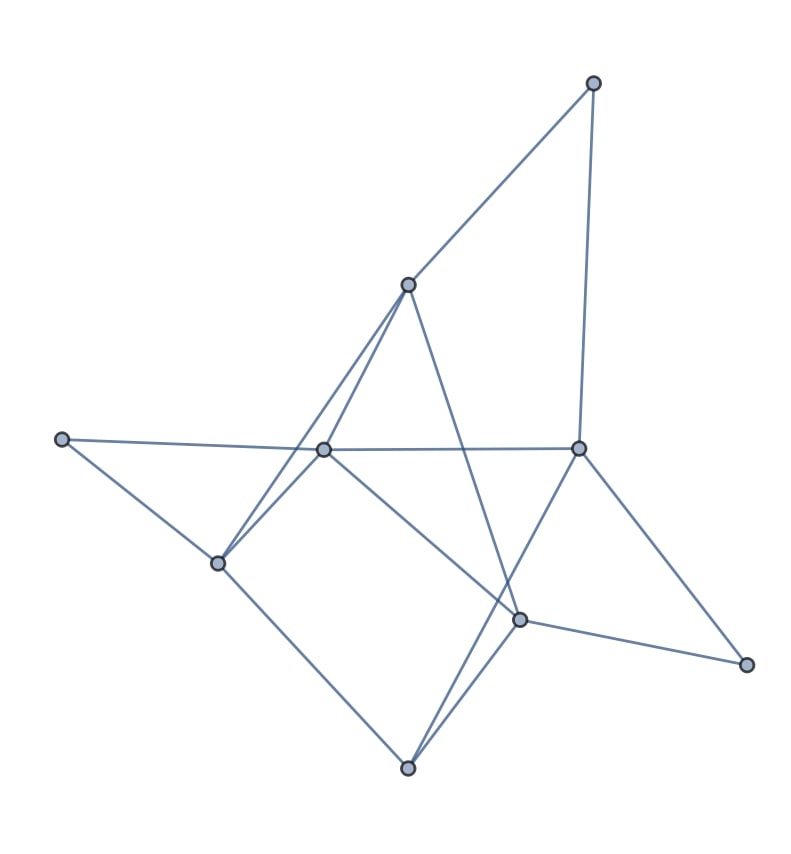}
    \caption{Graph with \texttt{graph6} string \texttt{HFRJIOY}.}
    \label{fig:HFRJIOY}
\end{figure}

\begin{figure}
    \centering
    \begin{minipage}{0.45\textwidth}
        \includegraphics[scale=0.15]{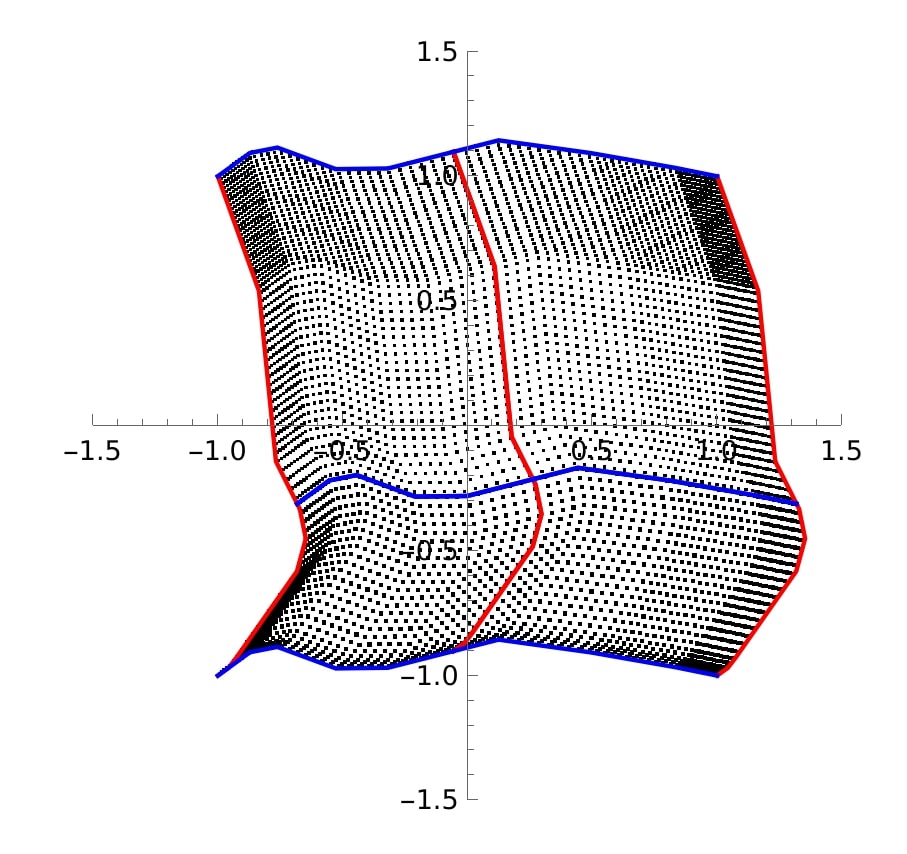} \subcaption*{$(-,+)$}
        \includegraphics[scale=0.15]{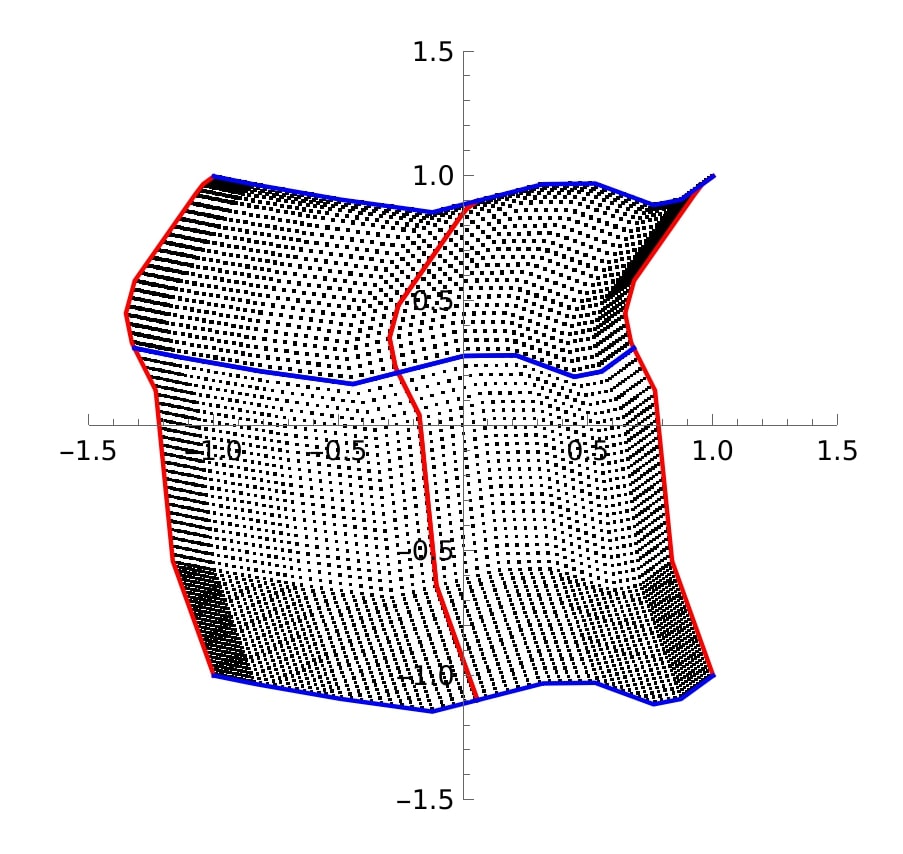} \subcaption*{$(-,-)$}
    \end{minipage}
    \hfil\hfil
    \begin{minipage}{0.45\textwidth}
        \includegraphics[scale=0.15]{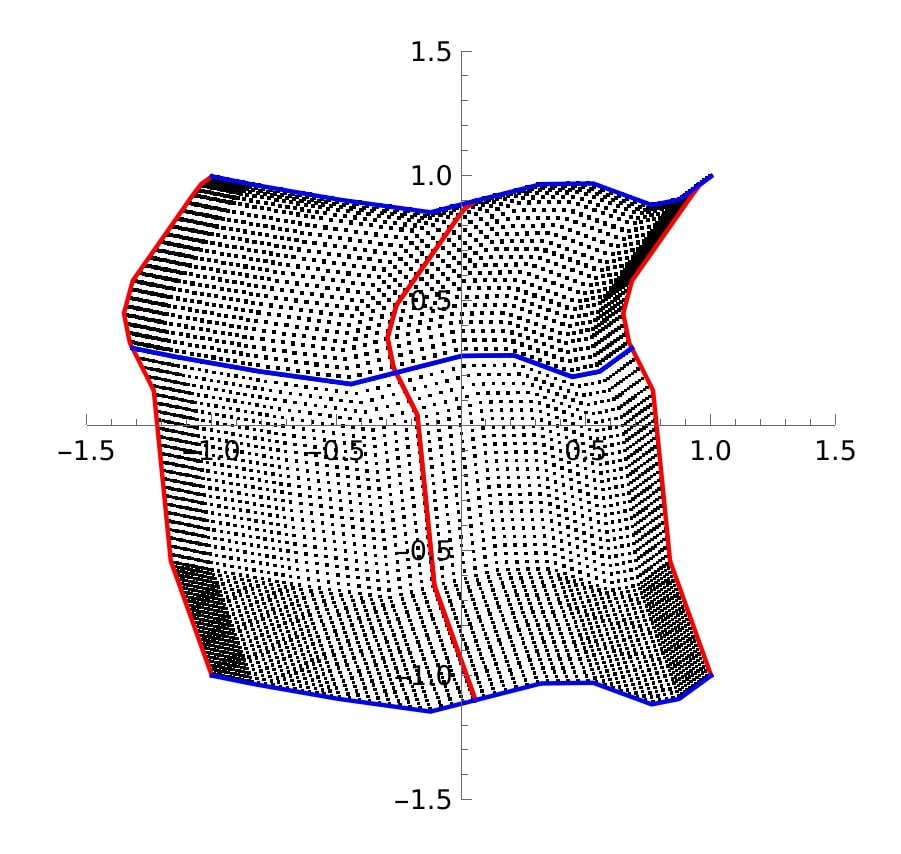} \subcaption*{$(+,+)$}
        \includegraphics[scale=0.15]{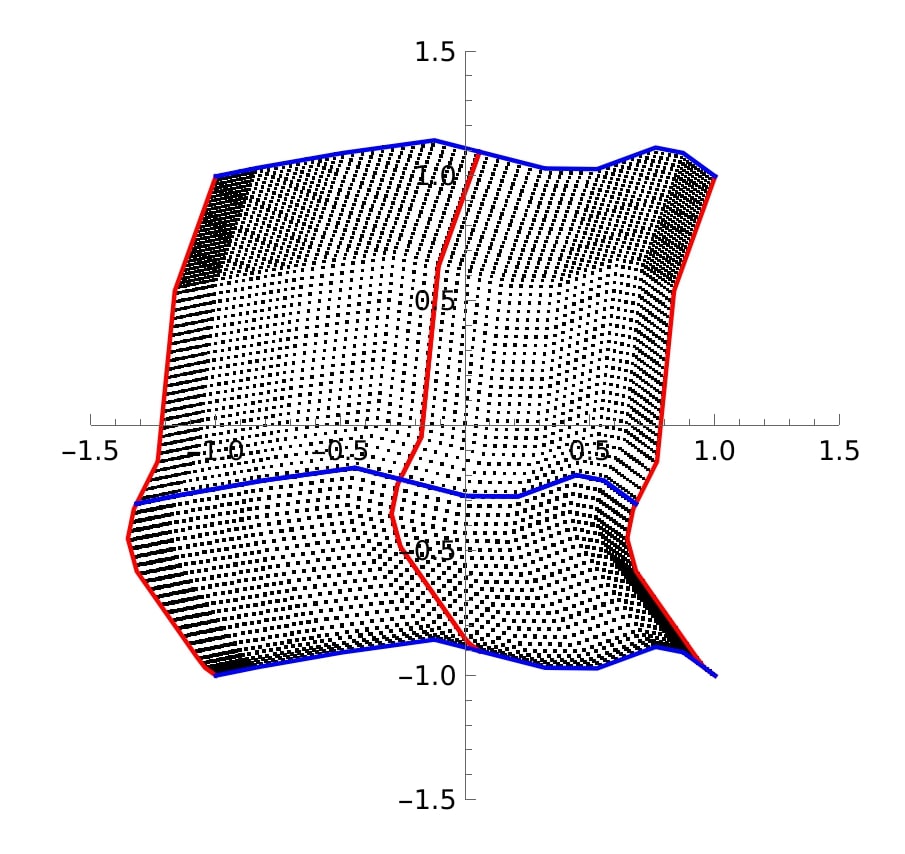} \subcaption*{$(+,-)$}
    \end{minipage}

    \caption{Images of the four possible orientations of the 2-dimensional 
Fiedler-carpet of the  graph in Figure~\ref{fig:HFRJIOY}.}
    \label{fig:fourientations}
\end{figure}

In Figure~\ref{fig:fourientations} we show the images of four different 
possible orientations of the Fiedler-carpet associated with the graph in 
Figure~\ref{fig:HFRJIOY}.
The $(+,+)$ orientation is shown in the top right, the $(-,-)$ in the bottom left, the $(+,-)$ in the bottom right, and the $(-,+)$ in the top left.
Notice that the ranges with the $(+,+)$ and the $(-,-)$ orientation are the 
same; the $(+,-)$ and $(-,+)$ orientations are reflections of each other 
across the vertical axis; 
the $(+,-)$ and $(+,+)$ orientations are reflections of each other 
across the horizontal axis; consequently the 
 $(-,+)$ and $(+,+)$ orientations are reflections of each other 
across the origin.


\begin{figure}
    \centering
    \includegraphics[scale=0.2]{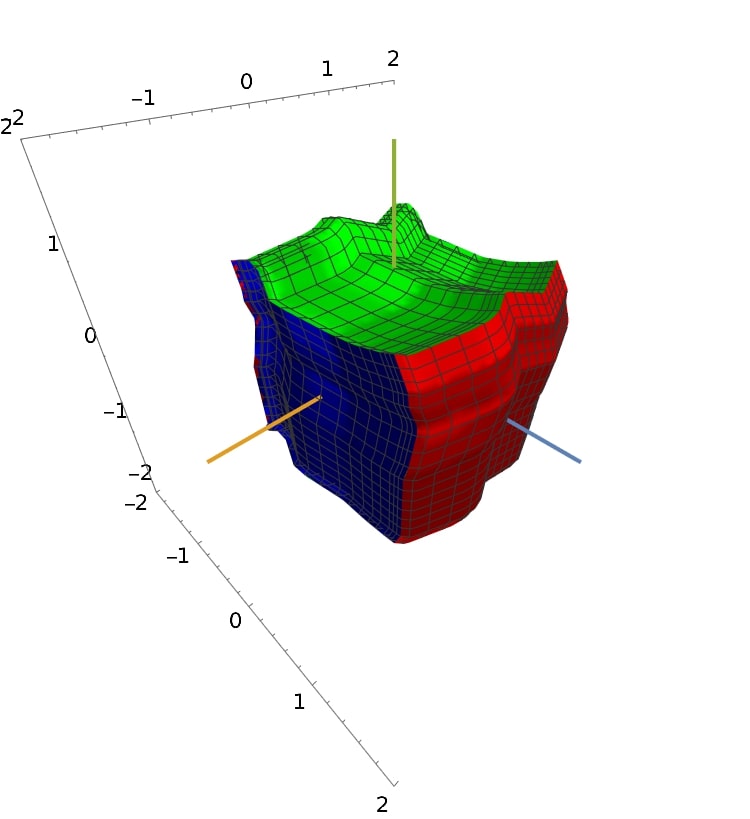}
    \includegraphics[scale=0.2]{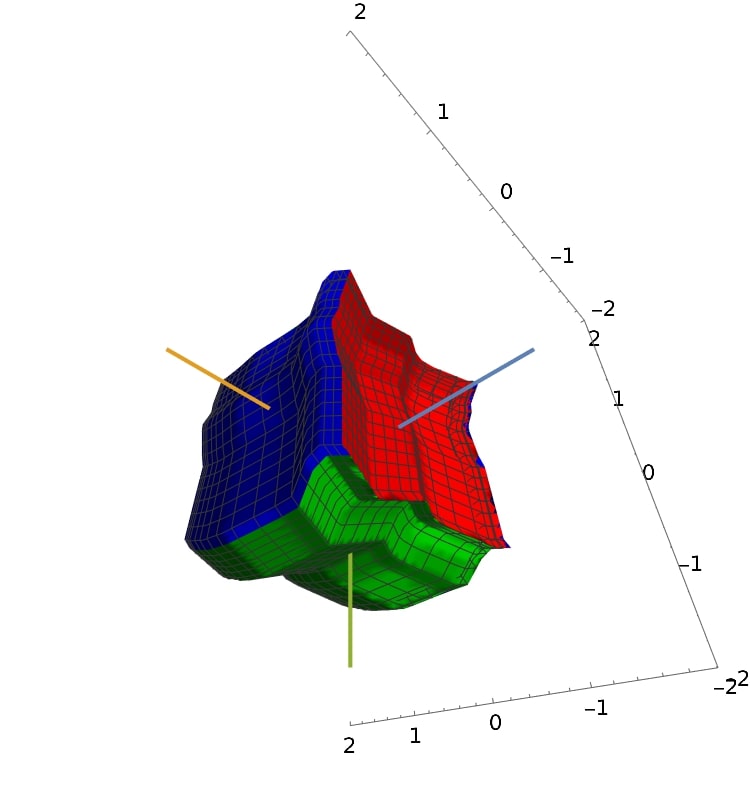}
    \caption{Image of the 3-dimensional Fiedler-carpet, 
in two different viewpoints, of the  graph in Figure~\ref{fig:HFRJIOY}.}
    \label{fig:3dimensions}
\end{figure}

We can also plot the map of the three-dimensional Fiedler-carpet 
of this graph.
In Figure~\ref{fig:3dimensions} different viewpoints are shown.
The blue, orange, and green lines are the coordinate axes in $\R^3$.
As can be seen from the images, they intersect at the origin \emph{inside} of the three-dimensional body.
Every plot we have created shows this intersection 
lying inside the map of the three-dimensional Fiedler-carpet.

If $k>2$, then $f:\R^k \to \R^k$ maps the $k$-dimensional hyper-rectangle with
vertices of $j$th coordinate $\min_i x_{ji}$ or $\max_i x_{ji}$ into the 
$k$-dimensional region with vertices of $j$th coordinate $\pm 1$.

Along the 
1-dimensional faces of this  $k$-dimensional range, 
all but one $a_i$ is fixed at its minimum/maximum. 
Without loss of generality, assume that 
$A<a_1 <B$.
Akin to the $k=2$ case, we are able to show that 
for each $A<a_1 <B$: 
$f_j (a_1 , M_{-1}) \ge 0$ and $f_j (a_1 , {\tilde M}_{-1}) \le 0$
$(j=2,\dots ,k)$,
where $M_{-1} = (\min_m x_{2m}, \dots , \min_m x_{km} )$  is the $(k-1)$-tuple  of 
the values of  
$a_2 ,\dots ,a_{k}$ all fixed at their minimum and
 ${\tilde M}_{-1} = (\max_m x_{2m}, \dots , \max_m x_{km})$  is the $(k-1)$-tuple  
of the values of 
$a_2 ,\dots ,a_{k}$ all fixed at their maximum values, respectively.

Indeed, for $j=2,\dots ,k$:
$$
\begin{aligned}
f_j (a_1 ,M_{-1} ) &=\sum_{i=1}^n d_i x_{ji} |x_{1i} -a_1 | +
    \sum_{l=2}^k  \sum_{i=1}^n d_i x_{ji} (x_{li} -\min_m  x_{lm} )  = \\
&=\sum_{i=1}^n d_i x_{ji} |x_{1i} -a_1 | +1 +  0 =\\ 
&=2\sum_{i\in H} d_i x_{ji} (x_{1i} -a_1 ) + 1 = \\
&=2\sum_{i\in {\bar H}} d_i x_{ji} (a_1 - x_{1i} ) + 1 
\end{aligned}
$$
as in the second double summation, only the term for $l=j$ is 1, the
others are zeros.  

Likewise,
$$
\begin{aligned}
f_j (a_1 ,{\tilde M}_{-1} ) &=\sum_{i=1}^n d_i x_{ji} |x_{1i} -a_1 | +
    \sum_{l=2}^k  \sum_{i=1}^n d_i x_{ji} (x_{li} -\max_m  x_{lm} )  =\\
 &=\sum_{i=1}^n d_i x_{ji} |x_{1i} -a_1 | -1 +  0 =\\ 
&=2\sum_{i\in H} d_i x_{ji} (x_{1i} -a_1 ) - 1 = \\
&=2\sum_{i\in {\bar H}} d_i x_{ji} (a_1 - x_{1i} ) - 1 .
\end{aligned}
$$
as in the second double summation, only the term for $l=j$ is -1, the
others are zeros.  
Note that counter-orienting $\uuu_1$ just results in $-a_1$ instead of
$a_1$ in the root of $f$.

The corresponding `corners'  are:
$f(A,M_{-1}) = (1,1,\dots ,1)$, 
$f(B,M_{-1}) = (-1,1,\dots ,1)$, 
$f(A,{\tilde M}_{-1}) = (1,-1,\dots ,-1)$, and
$f(B,{\tilde M}_{-1}) = (-1,-1,\dots ,-1)$.
These points are in one hyperplane, in which
the pieces of the parallel arcs
$f_j (a_1 , M_{-1})$ and $f_j (a_1 , {\tilde M}_{-1})$
sandwich the $f_1$ axis $(j=2,\dots ,k )$.
The same is true when another $a_j$ moves along a 1-dimensional face and the
others are fixed at their minima/maxima.
Therefore, the connected regions between these
parallel curves   sandwich the corresponding coordinate axes, respectively. 
Consequently, their intersection, which is subset of the whole 
connected region, contains the origin too. 

\section*{Appendix B}\label{appB}

\begin{figure}[htp]
	\begin{minipage}{.5\textwidth}
		\centering
		\includegraphics[width=6.8cm]{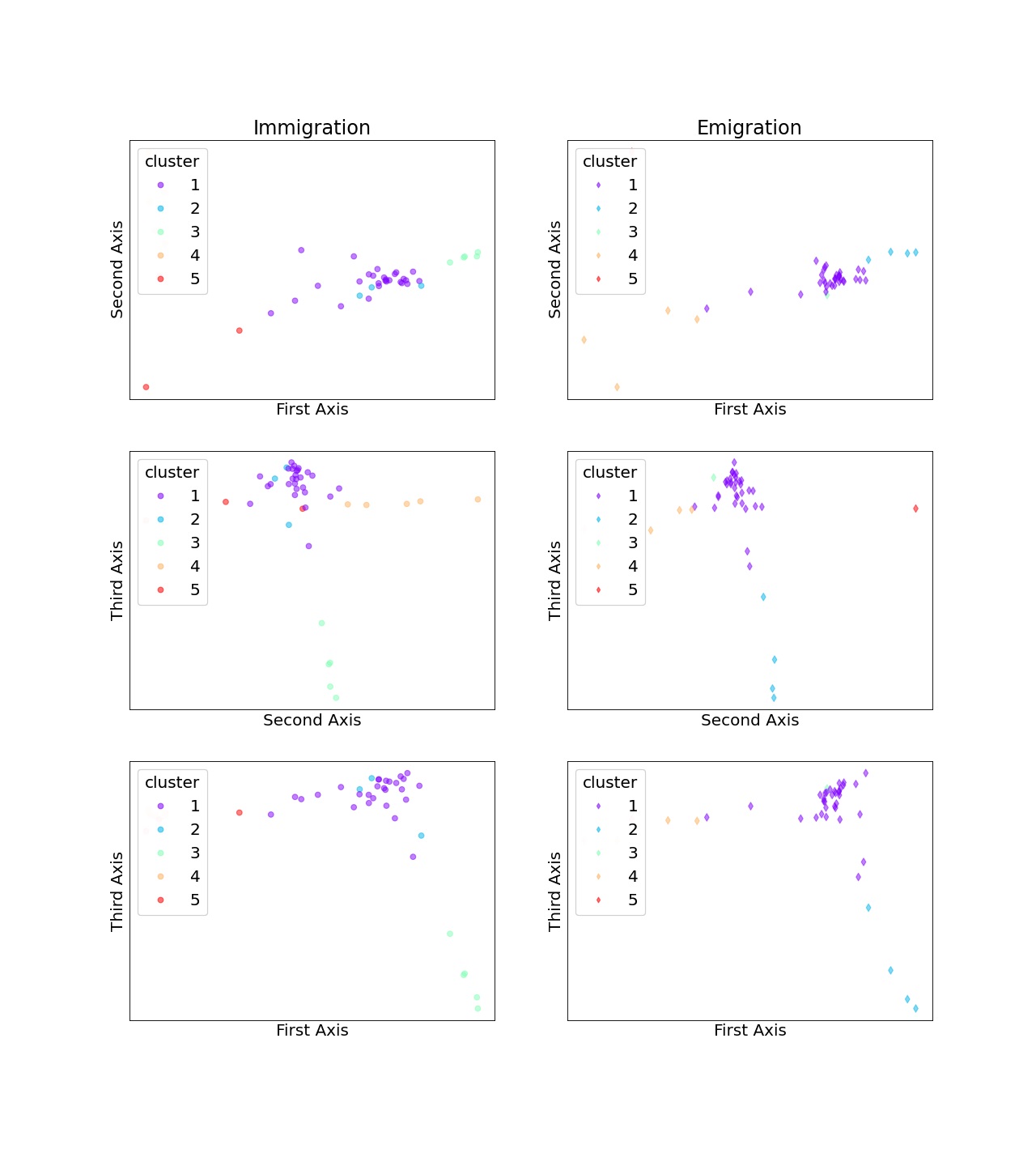}
		\subcaption{ 2015}
		\label{fig:2015_corresp}
	\end{minipage}		
	\hfil\hfil
	\begin{minipage}{.5\textwidth} 			
		\centering
		\includegraphics[width=6.8cm]{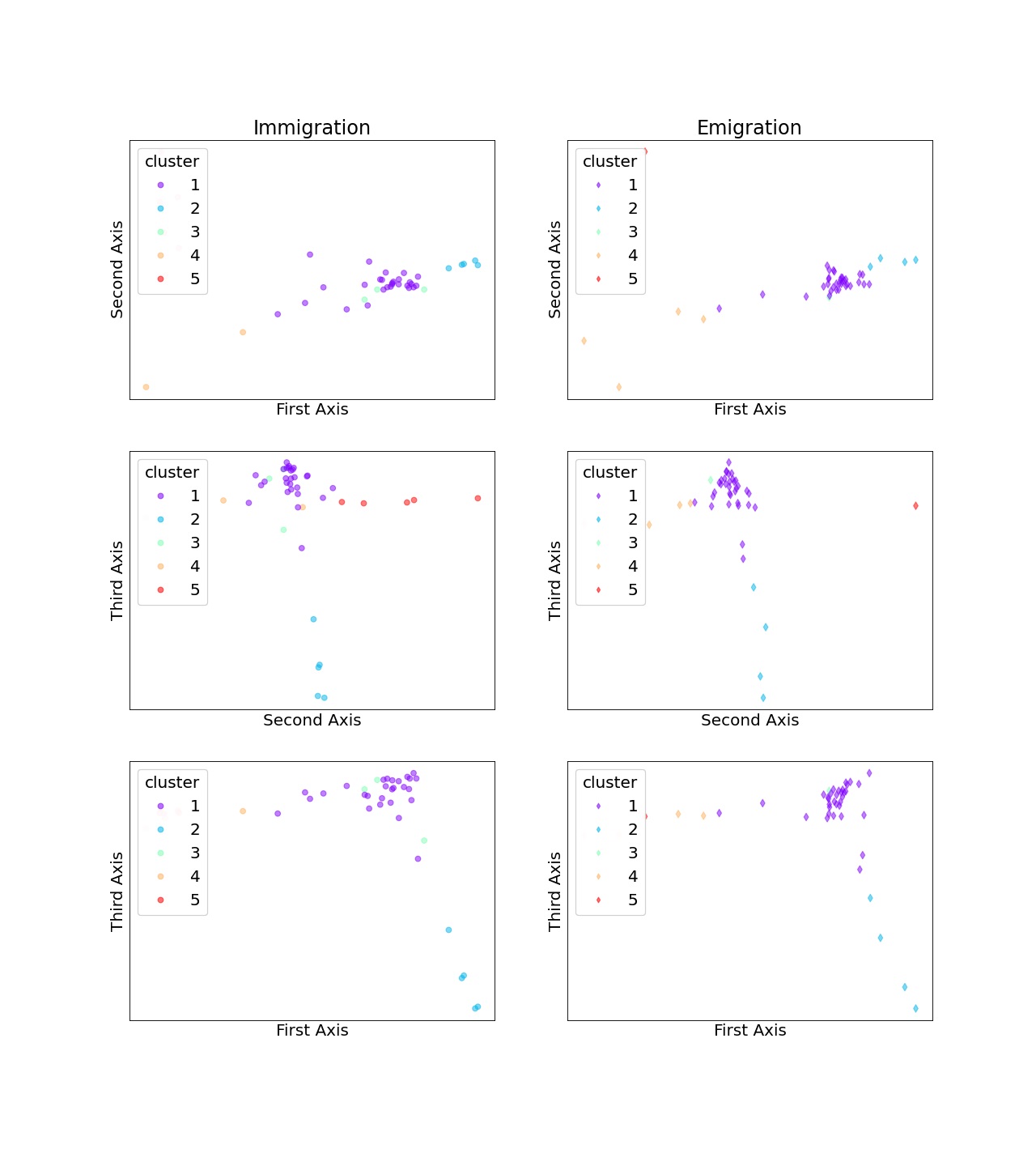}
		\subcaption{2019}
		\label{fig:2019_corresp}
	\end{minipage}			
\caption{Pairwise plots of the correspondence  analysis results based on the 
first three singular vector pairs, enhanced with the cluster memberships, 
illustrated by different colors.}
\end{figure}

Pairwise plots of the correspondence  analysis results based on the first three
coordinate axes (coordinates of the left singular vectors for imigration and 
right singular vectors for emigration data) are shown in 
Figures~\ref{fig:2015_corresp} and~\ref{fig:2019_corresp} for the two years.
The cluster memberships obtained in Section~\ref{appl} are
illustrated with different colors.

\end{document}